\documentclass[12pt]{article}
\usepackage{amsthm,amsmath,amssymb}
\def\8{\infty}
\def\sgn{\mbox{{\rm sgn}}}
\def\Res{\mathop{\hbox{\rm Res}}}

\newtheorem {thm}{Theorem}[section]
\newtheorem {lem}[thm]{Lemma}
\newtheorem {cor}[thm]{Corollary}

\newtheorem {prop}[thm]{Proposition}

\theoremstyle{remark}
\newtheorem{remark}[thm]{Remark}
\makeatletter

\@addtoreset{equation}{section}
\makeatother
\begin{document}
\title{A bilateral extension of the $q$-Selberg integral
}
\author{
M{\small ASAHIKO} I{\small TO}%
\footnote{
School of Science and Technology for Future Life,
Tokyo Denki University,
Tokyo 120-8551, Japan 
\quad email: {\tt mito@cck.dendai.ac.jp}}
\ and P{\small ETER} J. F{\small ORRESTER}%
\footnote{
Department of Mathematics and Statistics, The University of Melbourne, Victoria 3010, Australia 
\quad email: {\tt p.forrester@ms.unimelb.edu.au}
}
}
\date{}
\maketitle
\begin{abstract}
A multi-dimensional bilateral $q$-series extending the $q$-Selberg integral is studied 
using concepts of {\it truncation}, {\it regularization} and {\it connection formulae}.
Following Aomoto's method, which involves regarding the $q$-series as a solution of a $q$-difference equation fixed by its asymptotic behavior, 
an infinite product evaluation is obtained. 
The $q$-difference equation is derived 
applying the {\it shifted symmetric polynomials} introduced by Knop and Sahi. 
As a special case of the infinite product formula, Askey--Evans's $q$-Selberg integral evaluation 
and its generalization by 
Tarasov--Varchenko and Stokman is reclaimed, and an explanation in the context of Aomoto's setting is thus provided.
\end{abstract}

{\scriptsize 2010 {\it Mathematical Subject Classification.} Primary 33D15, 33D67; Secondary 39A13}

{\scriptsize  {\it Keywords.} 
Askey--Evans's Selberg integral;
Aomoto's method; 
connection formulae; 
Knop--Sahi's shifted symmetric polynomials
}
\section{Introduction}
The Selberg integral  \cite{Se44} is a multi-dimensional generalization of the evaluation of the Euler beta integral in terms of products of gamma functions. It reads
\begin{eqnarray}
\label{eq:S}
&&\frac{1}{n!}\int_{0}^{1}\cdots\int_{0}^{1}
\prod_{i=1}^n z_i^{\alpha-1}(1-z_i)^{\beta-1}\prod_{1\le j<k\le n}|z_j-z_k|^{2\tau}\,dz_1dz_2\cdots dz_n
\nonumber\\
&&\quad =\prod_{j=1}^n\frac{\Gamma(\alpha+(j-1)\tau)\Gamma(\beta+(j-1)\tau)\Gamma(j\tau)}{\Gamma(\alpha+\beta+(n+j-2)\tau)\Gamma(\tau)}.
\end{eqnarray}
The Euler beta integral is fundamental to the theory of hypergeometric functions.
And  fundamental in the theory of hypergeometric functions is the notion of a $q$ (or basic)
generalization.
From the 1980's to 1990's,  a $q$-analog of the Selberg integral formula was established and proved 
by Askey \cite{As80}, Habsieger \cite{Ha88}, Kadell \cite{Kad88}  and Evans \cite{Ev92}. 
\begin{prop}
\label{prop:AHKE}
Suppose $|q^\alpha|<1$, $|q^\beta|<1$ and $\tau$ is a positive integer. Let the $q$-integral
be defined in terms of a sum as specified by {\rm (\ref{2.4a})} below. We have
\begin{eqnarray}
\label{eq:AHKE}
&&\frac{1}{n!}\int_0^1\cdots\int_0^1\prod_{i=1}^n
z_i^{\alpha-1}(qz_i)_{\beta-1}
\prod_{1\le j<k\le n\atop 1-\tau\le l\le \tau-1}
(z_j-q^lz_k)
\prod_{1\le j<k\le n}
(z_j-z_k)
%
\,d_qz_1\cdots d_qz_n\nonumber\\
&&\qquad=q^{\alpha\tau{n\choose 2}+2\tau^2{n\choose 3}}
\prod_{j=1}^n \frac
{\Gamma_q(\alpha+(j-1)\tau)\Gamma_q(\beta+(j-1)\tau)\Gamma_q(j\tau)}
{\Gamma_q(\alpha+\beta+(n+j-2)\tau)\Gamma_q(\tau)}.
\end{eqnarray}
\end{prop}
An important point for present purposes is that the above $q$-generalization of the Selberg
integral is restricted to the case that the parameter $\tau$ is a positive integer.

On the other hand, in the 1990's, Aomoto gave a 
$q$-analogue of (\ref{eq:S}) valid for 
general complex  $\tau$.
\begin{prop}
\label{prop:aomoto} 
{\rm \cite[p.\,121,\,Proposition 2]{Ao98}}
Let $\alpha,\beta$ and $\tau$ be complex numbers satisfying \\
$
|q^{\alpha+(i-1)\tau}|<1$ for 
$i=1,2,\ldots,n$. 
Then
\begin{eqnarray}
\label{eq:I(zeta)}
&&\hskip -19pt
\int_{z_1=0}^{1}
\int_{z_2=0}^{q^{\tau}z_1}
\!\!\cdots\!\!
\int_{z_n=0}^{q^{\tau}z_{n-1}}
\prod_{i=1}^nz_i^\alpha\frac{(qz_i)_\8}{(q^\beta z_i)_\8}
\prod_{1\le j<k\le n}
z_j^{2\tau-1}\frac{(q^{1-\tau}z_k/z_j)_\8}{(q^\tau z_k/z_j)_\8}(z_j-z_k)
\frac{d_qz_n}{z_n}\cdots 
\frac{d_qz_2}{z_2}
\frac{d_qz_1}{z_1}
\nonumber\\
&&\qquad=q^{\alpha\tau{n\choose 2} +2\tau^2{n\choose 3}}\prod_{j=1}^n 
\frac
{\Gamma_q(\alpha+(j-1)\tau)\Gamma_q(\beta+(j-1)\tau)\Gamma_q(j\tau)}
{ 
\Gamma_q(\alpha +\beta+(n+j-2)\tau)\Gamma_q(\tau)
}.
\end{eqnarray}
\end{prop}

An immediate question presents itself: how are (\ref{eq:AHKE}) and (\ref{eq:I(zeta)}) 
related when in the latter $\tau$ is a positive integer? First, one should remark that in this case
they both 
coincide with Selberg's formula (\ref{eq:S}) in the limit $q\to 1$.
However, this fact aside, it soon becomes apparent that it is not possible to obtain one result 
from the other by analytic continuation. How to remedy this situation motivates the first study of this paper. Thus we present in Section \ref{section:01} the
more general theory of Aomoto relating to (\ref{eq:I(zeta)}), which allows for the development
of a theory of the {\it truncated Jackson integral of type A}. This theory is based on $q$-difference
equations and asymptotic behaviors.
In particular, by adding our supplementary results (see Lemmas \ref{lem:I-vanishing} and \ref{lem:I(zeta)=I()}), we are able to obtain Proposition \ref{prop:AHKE} as a corollary of Proposition \ref{prop:aomoto}.

However our main purpose of this paper is not just a supplementary commentary to Aomoto's work. 
In the paper \cite{As80}, Askey gave a conjecture for another $q$-extension of the Selberg integral, 
which is proved by Evans \cite{Ev94} using Anderson's method \cite{An91}. 
Askey--Evans's formula is written as follows:
\begin{prop}
\label{prop:AE}
{\rm \cite[p.\,342, Theorem 1, (1.9)]{Ev94}} If $\tau$ is a positive integer, then
\begin{eqnarray}
\label{eq:Evans0}
&&\frac{1}{n!}\int_{x_1}^{x_2}\cdots\int_{x_1}^{x_2}
\prod_{i=1}^n\Big(\frac{qz_i}{x_1}\Big)_{\!\!\alpha-1}
\Big(\frac{qz_i}{x_2}\Big)_{\!\!\beta-1}
\!\!
\prod_{1\le j<k\le n\atop 1-\tau\le l\le \tau-1}
(z_j-q^lz_k)
\prod_{1\le j<k\le n}
(z_j-z_k)
\,d_qz_1\cdots d_qz_n\nonumber\\
\\[-14pt]
&&=(-1)^{\tau {n \choose 2}}q^{\tau^2{n \choose 3}-{\tau \choose 2}{n \choose 2}}
\nonumber\\
&&\ \times
\prod_{j=1}^n\frac{\Gamma_q(\alpha+(j-1)\tau)\Gamma_q(\beta+(j-1)\tau)\Gamma_q(j\tau)}
{\Gamma_q(\alpha+\beta+(n+j-2)\tau)\Gamma_q(\tau)}
\frac{(x_1x_2)^{1+(j-1)\tau}
}{x_1-x_2}
\Big(\frac{x_2}{x_1}\Big)_{\!\!\alpha+(j-1)\tau}\Big(\frac{x_1}{x_2}\Big)_{\!\!\beta+(j-1)\tau}.
\nonumber
\end{eqnarray}
\end{prop}
This is the $q$-analog of the Selberg integral  (\ref{eq:S})  whose integral area is transformed from $[0,1]^n$ to $[x_1,x_2]^n$
by using the linear transformation $z_i\mapsto (x_2-x_1)z_i+x_1$. 
This formula, like  (\ref{eq:AHKE}),
was also formulated under the assumption the parameter $\tau$ is a positive integer. 
Tarasov and Varchenko \cite{TV97} and Stokman \cite{St00} independently 
gave an extension of Askey--Evans's formula in the case of 
$\tau$ being an arbitrary complex number (see (\ref{eq:TV,S}) in Corollary \ref{cor:main1}), 
using a residue calculus on a certain contour integral 
(see also Gustafson's $q$-Selberg contour integral \cite{Gu90,Gu94}).

Our primary goal is to investigate the case of $\tau$ being an arbitrary complex number
defining a bilateral extension of the $q$-Selberg integral,
following Aomoto's method as presented in Section \ref{section:01}, i.e., 
the $q$-difference equations and its solutions fixed by the asymptotic behaviors.
The main theorem of this paper, providing the solution of this problem,
 is Theorem \ref{thm:main1}. 
(Theorem \ref{thm:main1} is equivalent to Theorem \ref{thm:main2}, 
whose expression seems to be simpler by the term {\em regularization}.)
And as with our findings in Section \ref{section:01}
we can understand 
the formula (\ref{eq:TV,S}) of Tarasov--Varchenko and Stokman
as a special case 
(called the {\it truncation}) 
of the formula (\ref{eq:main1}) in Theorem \ref{thm:main1}, 
and can also understand Askey--Evans's formula (\ref{eq:Evans0}) 
as a special case of 
(\ref{eq:Evans3}) in Corollary \ref{cor:main1+}, i.e., 
the formula (\ref{eq:main1}) with the restriction of $\tau$ being a positive integer. (See Section \ref{subsection:02.1}.)
In particular, we can see the degeneration occurs in the factors written by theta functions in the left-hand side of the formula (\ref{eq:main1}). 
This means it is not so easy to guess the exact form of the formula (\ref{eq:main1}) only from knowledge of Askey--Evans's formula (\ref{eq:Evans0}) as a special case, rather the Aomoto
viewpoint plays an essential role. \\ 

The method for proving the results in this paper is consistent with the concept introduced by Aomoto 
and Aomoto--Kato in the early 1990's
in the series of papers \cite{Ao90,Ao91,Ao94,Ao95-1,Ao98,AK91,AK93,AK94-1,AK94-2}. 
Aomoto showed an isomorphism between a class of Jackson integrals of hypergeometric type,
which he called the {\it $q$-analog de Rham cohomology} \cite{Ao90,Ao91}, and a class of theta functions, i.e., 
holomorphic functions possessing a quasi-periodicity \cite[Theorem 1]{Ao95-1}. 
This isomorphism indicates that it is essential to analyze the class of holomorphic functions 
as a counterpart of that of Jackson integrals 
in order to know the structure of $q$-hypergeometric functions, in particular, 
the meaning of known special formulae. 
In this paper the process to obtain the holomorphic functions through this isomorphism is called the {\it regularization}. 
When we fix a basis of the class of holomorphic functions as a linear space, 
an arbitrary function of the space can be expressed as a linear combination of the elements of the specific basis, 
which he called the {\it connection formula} \cite[Theorem 3]{Ao95-1}. 
As its simplest examples, Ramanujan's $_1\psi_1$ summation formula and 
the $q$-Selberg integral \cite{As80,Ha88,Kad88,Ev92}  have been explained. See the original literature \cite[Examples 1, 2]{Ao95-1} and the recent review \cite{IF13} for details.
One way to choose a good basis is to fix it by its asymptotic behavior of a limiting process with  
respect to parameters included in the definition of the Jackson integral of hypergeometric type.
And the asymptotic behavior can be calculated from the Jackson integrals possessing 
appropriate cycles which include their critical points. We call the process to fix the cycles {\it truncation}.
(These cycles are called the {\it characteristic cycles} \cite{AK94-2} 
or the {\it $\alpha$-stable or $\alpha$-unstable cycles} \cite{Ao94} by Aomoto. 
The meaning of ``$\alpha$" is mentioned in Section \ref{section:01}. 
The word {\it truncation} itself is first used by van Diejen in other context \cite{vD97, Ito06-2}.)
The connection formula is also characterized as a formula showing that 
a multi-dimensional bilateral series giving a general solution 
of the $q$-difference equation of the Jackson integrals with respect to parameters 
is expressed as a linear combination of multi-dimensional unilateral series
as special solutions, each fixed by their asymptotic behaviors \cite[Theorem (4.2)]{Ao94}. 
(We can see different examples of $q$-difference equations and the connection formulae in \cite{Ito08,Ito09,IN15,IS08}, 
and \cite{IN15,IS08} explain the Sears--Slater transformation for the very-well-poised $q$-hypergeometric series 
from the present view point in the setting of $BC$ type symmetry.)

The paper is organized as follows. 
After defining some basic terminology in Section \ref{section:00}, 
we first show the product expression of the $q$-Selberg integral along Aomoto's setting 
(we called it the {\em Jackson integral of A-type}) 
using concepts of {\it truncation}, {\it regularization} and {\it connection formulae} in Section \ref{section:01}. 
Though the Jackson integral of A-type can be obtained from our other example, 
we explain it individually, because the Jackson integral of A-type has simpler structure than the other Jackson integral,
it is instructive in outlining the concept of this paper, and it highlights the issue of the relationship
between (\ref{eq:AHKE}) and (\ref{eq:I(zeta)}).
Section \ref{section:02} is devoted to explaining a bilateral extension of Askey--Evans's $q$-Selberg integral.
Its situation looks a little more complex  
than the case of the Jackson integral of A-type in their details, but still it is consistent with 
the outlines of the proofs for the product expressions of these sums. 
In the Appendix we explain the detail of the derivation of the $q$-difference equation. 
In particular, we applied the {\it shifted symmetric polynomials} introduced by Knop and Sahi \cite{KS96} 
to the key lemma (Lemma \ref{lem:Se=Se}) for deriving the $q$-difference equation.

\section{Definition of the Jackson integral}
\label{section:00}
Throughout this paper, we fix $q$ as $0<q<1$ and use the symbols 
$(a)_\8:=\prod_{i=0}^\8(1-q^i a)$ and $(a)_N:=(a)_\8/(q^Na)_\8$. 
We define $\theta(a)$ by $\theta(a):=(a)_\8(q/a)_\8$, 
which satisfies 
\begin{equation}
\label{eq:00quasi-period}
\theta(qa)=-\theta(a)/a.
\end{equation}
By repeated use of the latter, $\theta(a)$ satisfies 
\begin{equation}
\label{eq:00quasi-period2}
\theta(a)/\theta(q^m a)=(-a)^mq^{{m\choose 2}}
\quad\mbox{for}\quad m\in \mathbb{Z}.
\end{equation}
We define $\Gamma_q(x)$ by $\Gamma_q(x):=(1-q)^{1-x}(q)_\infty/(q^x)_\infty$, 
which satisfies
$$
\Gamma_q(x)\Gamma_q(1-x)=\frac{(1-q)(q)_\infty^2}{(q^x)_\infty(q^{1-x})_\infty}
=(1-q)\frac{(q)_\infty^2}{\theta(q^x)},
$$
this being a $q$-analog of the relation 
$
\Gamma(x)\Gamma(1-x)=\pi/\sin \pi x$. \\

Let $S_n$ be the symmetric group on $\{1,2,\ldots, n\}$.
For a function $f(z)=f(z_1,z_2,\ldots,z_n)$ on $({\mathbb{C}^*})^n$, 
we define an action of the symmetric group $S_n$ on $f(z)$ by 
$$
(\sigma f)(z):=f(\sigma^{-1}(z))=f(z_{\sigma(1)},z_{\sigma(2)},\ldots,z_{\sigma(n)})
\quad\mbox{for}\quad \sigma\in S_n.
$$
We say that a function $f(z)$ on $({\mathbb{C}^*})^n$
is {\it symmetric} or {\it skew-symmetric} 
if $\sigma f(z)=f(z)$ or $\sigma f(z)=(\sgn\,\sigma )\,f(z)$ 
for all $\sigma \in S_n$, respectively.
We denote by ${\cal A} f(z)$ 
the alternating sum over $S_n$ defined by 
\begin{equation}
\label{eq:00Af}
({\cal A} f)(z):=\sum_{\sigma\in S_n}(\sgn\, \sigma)\,\sigma f(z),
\end{equation}
which is skew-symmetric. \\
\par
For $a,b\in \mathbb{C}$, we define 
\begin{equation}
\label{eq:00jac1}
\int_a^b f(z)d_qz:=\int_0^b f(z)d_qz-\int_0^a f(z)d_qz,
\end{equation}
where
\begin{equation}\label{2.4a}
\int_0^a f(z)d_qz:=(1-q)\sum_{\nu=0}^\8 f(aq^\nu)aq^\nu,
\end{equation}
which is called the {\it Jackson integral}. As $q\to 1$, 
$\int_a^b f(z)d_qz\to \int_a^b f(z)dz$ \cite{AAR99}. 
In this paper we use the Jackson integral of multiplicative measure, specified by
$$
\int_0^a f(z)\frac{d_qz}{z}=(1-q)\sum_{\nu=0}^\8 f(aq^\nu). 
$$
as is consistent with (\ref{2.4a}). 
Let $\mathbb{N}$ be the set of non-negative integers. 
For a function $f(z)=f(z_1,\ldots,z_n)$ on $(\mathbb{C}^*)^n$ and 
an arbitrary point $x=(x_1,\ldots,x_n)\in (\mathbb{C}^*)^n$, 
we define the multiple Jackson integral as
\begin{equation}
\label{eq:00jac2}
\int_0^{\mbox{\small $x$}}f(z)\,\frac{d_qz_1}{z_1}\wedge\cdots\wedge\frac{d_qz_n}{z_n}
:=(1-q)^n\sum_{(\nu_1,\ldots,\nu_n)\in \mathbb{N}^n}f(x_1 q^{\nu_1},\ldots,x_n q^{\nu_n}). 
\end{equation}
In this paper we use the bilateral sum extending the Jackson integral (\ref{eq:00jac1})
\begin{equation}
\label{eq:00jac3}
\int_{a\infty}^{b\infty} f(z)d_qz:=\int_0^{b\infty} f(z)d_qz-\int_0^{a\infty} f(z)d_qz,
\end{equation}
where
$$
\int_0^{a\infty} f(z)d_qz:=(1-q)\sum_{\nu=-\infty}^\infty f(aq^\nu)aq^\nu,
\ \mbox{ i.e., }\ 
\int_0^{a\infty} f(z)\frac{d_qz}{z}:=(1-q)\sum_{\nu=-\infty}^\infty f(aq^\nu).
$$
We also use the multiple bilateral sum extending the Jackson integral (\ref{eq:00jac2})
\begin{equation}
\label{eq:00jac4}
\int_0^{\mbox{\small $x$}\8}f(z)\,\frac{d_qz_1}{z_1}\wedge\cdots\wedge\frac{d_qz_n}{z_n}
:=(1-q)^n\sum_{(\nu_1,\ldots,\nu_n)\in \mathbb{Z}^n}f(x_1 q^{\nu_1},\ldots,x_n q^{\nu_n}),
\end{equation}
which we also call the {\it Jackson integral}. By definition the Jackson integral (\ref{eq:00jac4}) 
is invariant under the shift $x_i\to qx_i$, $1\le i\le n$.
While we can consider the limit $q\to1$ for the Jackson integral (\ref{eq:00jac2}) defined over $\mathbb{N}^n$,  
the Jackson integral (\ref{eq:00jac4}) defined over $\mathbb{Z}^n$ generally diverges if $q\to1$. 
However, as we will see later, since the {\it truncation} of 
the Jackson integral (\ref{eq:00jac4}) is corresponding to the sum (\ref{eq:00jac2}) over $\mathbb{N}^n$,  
if we need to consider the limit $q\to1$, we switch from (\ref{eq:00jac4}) to (\ref{eq:00jac2}) by the process of the truncation. 
For simplicity of notation, we use the symbol 
$$\varpi_q
=\frac{d_qz_1}{z_1}\wedge\cdots\wedge\frac{d_qz_n}{z_n}.$$

\section{Aomoto's $q$-extension of the Selberg integral}
\label{section:01}
In this section we will review some known results in the context of Aomoto's $q$-extension of the Selberg integral. 

\subsection{Aomoto's setting}
For $\alpha\in \mathbb{C}$, $a_1,b_1,t\in \mathbb{C}^*$, $z=(z_1,z_2,\ldots,z_n)\in (\mathbb{C}^*)^n$,  
let $\Phi(z)$ and $\Delta(z)$ be specified by 
\begin{eqnarray}
&&\Phi(z):=
\prod_{i=1}^nz_i^\alpha\frac{(qa_1^{-1}z_i)_\8}{(b_1 z_i)_\8}
\prod_{1\le j<k\le n}
z_j^{2\tau-1}\frac{(qt^{-1}z_k/z_j)_\8}{(t z_k/z_j)_\8},
\label{eq:Phi0}\\
&&\Delta(z):=\prod_{1\le i<j\le n}(z_i-z_j), 
\label{eq:Delta}
\end{eqnarray}
where $\tau$ is given by $t=q^{\tau}$. 
For $x=(x_1,x_2\ldots,x_n)\in (\mathbb{C}^*)^n$ we define $I(x)$ by
\begin{equation}
\label{eq:I(xi)}
I(x)=I(x_1,x_2,\ldots,x_n):=\int_0^{\mbox{\small $x$}\8}
\Phi(z)
\Delta(z)
\varpi_q,
\end{equation}
which is called the {\em Jackson integral of {\rm A}-type} in the context of \cite{Ao98}. 
For a general point $x\in  (\mathbb{C}^*)^n$, excluding poles of $\Phi(z)$, 
$I(x)$ converges absolutely under the condition 
\begin{equation}
\label{eq:convergence1-I}
|qa_1^{-1}b_1^{-1}|<|q^{\alpha}t^{2i-2}|<1\quad\mbox{for}\quad i=1,2,\ldots,n.
\end{equation}
Let $\zeta$ be the point defined by 
\begin{equation}
\label{3.4'}
\zeta:=(a_1,a_1t,a_1t^2,\ldots,a_1t^{n-1})\in (\mathbb{C}^*)^n,
\end{equation}
and let $\Lambda$ be the subset of $\mathbb{Z}^n$ defined by \\[10pt]
\hskip 95pt$
\Lambda:=\{(\nu_1,\nu_2,\ldots,\nu_n)\in \mathbb{Z}^n\,;\, 0\le\nu_1\le\nu_2\le \cdots\le \nu_n\}.
$\\[-5pt]
The set $\Lambda$ is written as $\Lambda=\{\sum_{i=1}^n m_i\epsilon_i\,;\, m_i\in \mathbb{N}\}\cong\mathbb{N}^n$
where $\epsilon_i=(\hskip 1pt\overbrace{0,\ldots,0}^{i-1},1,\ldots,1)\in \mathbb{Z}^n$. 
Since $\Phi(\zeta q^{\nu})=\Phi(a_1q^{\nu_1},a_1tq^{\nu_2},\ldots,a_1t^{n-1}q^{\nu_n})=0$ 
if $\nu\not\in \Lambda$, 
by definition $I(\zeta)$ is defined as the sum over the fan region $\Lambda\cong\mathbb{N}^n$. 
For this special point $\zeta$, $I(\zeta)$ converges absolutely if 
\begin{equation}
\label{eq:convergence2-I}
|q^{\alpha}t^{i-1}|<1\quad\mbox{for}\quad i=1,2,\ldots,n
\end{equation}
and  we call $I(\zeta)$ the {\em truncated Jackson integral of {\rm A}-type}.
Note that $I(\zeta)$ is expressed as the following iterated $q$-integral form:
$$
I(\zeta)=
\int_{z_1=0}^{a_1}
\int_{z_2=0}^{tz_1}
\cdots
\int_{z_n=0}^{tz_{n-1}}
\Phi(z)
\Delta(z)
\frac{d_qz_n}{z_n}\cdots 
\frac{d_qz_2}{z_2}
\frac{d_qz_1}{z_1}.
$$

Although we regard all parameters as complex numbers throughout the paper, 
it is often very important to distinguish between 
the parameter $\tau$ being a positive integer or not, as we see 
in Propositions \ref{prop:AHKE} and \ref{prop:AE}  for instance. 
We initially state a basic property of $I(x)$ under the assumption $\tau$ is not a positive integer.
Let $\mathbb{Z}_+$ be the set of positive integers. \\

By definition the function $\Phi(z)$ 
satisfies the {\em quasi-symmetric property} that 
$$
\sigma \Phi (z)=U_{\sigma}(z)\Phi(z) \quad\mbox{for}\quad \sigma \in S_n, 
$$
where
\begin{equation}
\label{eq:Usigma}
U_\sigma(z):=\prod_{1\le i<j\le n\atop \sigma^{-1}(i)>\sigma^{-1}(j)} 
\Big(\frac{z_i}{z_j}\Big)^{\!\!1-2\tau}\frac{\theta(q^{1-\tau}z_i/z_j)}{\theta(q^\tau z_i/z_j)},
\end{equation}
which is invariant under the $q$-shift $z_i\to q z_i$. 
From (\ref{eq:Usigma}) and $\sigma\Delta(z)=(\sgn\,\sigma)\Delta(z)$, we have
\begin{equation}
\label{eq:Usigma-I}
\sigma I(x)=(\sgn\,\sigma)U_\sigma(x)I(x).
\end{equation}

\begin{lem}
\label{lem:I-vanishing}
Suppose $\tau\not\in\mathbb{Z}_+$. 
If $x_i=x_j$ for some $i$ and $j$ $(1\le i<j\le n)$, then $I(x_1,x_2,\ldots,x_n)=0$. 
\end{lem}
\noindent
{\bf Proof.} 
Set $\sigma$ as the interchange of $i$ and $j$. If we impose $x_i=x_j$, then $\sigma I(x)=I(x)$,
so that we have  
$(1+ U_\sigma(x))I(x)=0$ 
from (\ref{eq:Usigma-I}). Since $(1+ U_\sigma(x))\ne 0$, we obtain $I(x)=0$. $\square$\\

On the other hand, 
under the assumption that $\tau$ is a positive integer 
we generally have $I(x_1,\ldots,x_n)\ne 0$ even if $x_i=x_j$ $(1\le i<j\le n)$. 
In particular, we have the following:

\begin{lem}
\label{lem:I(zeta)=I()}
Suppose that $\tau\in \mathbb{Z}_+$. For an arbitrary $x\in \mathbb{C}^*$
\begin{equation}
\label{eq:I(zeta)=I()}
I(x,xt,\ldots,xt^{n-1})=\frac{I(x,x,\ldots,x)}{n!}. 
\end{equation}
\end{lem}
\noindent
{\bf Proof.} Under the condition $\tau\in \mathbb{Z}_+$, since  
\begin{eqnarray}
&&\prod_{1\le j<k\le n}z_j^{2\tau-1}\frac{(q^{1-\tau}z_k/z_j)_\infty}{(q^\tau z_k/z_j)_\infty}(z_j-z_k)
=
\prod_{1\le j<k\le n\atop 1-\tau\le l\le \tau-1}
(z_j-q^lz_k)
\prod_{1\le j<k\le n}
(z_j-z_k)
\nonumber\\
&&\qquad=
(-1)^{\tau{n\choose 2}}q^{-{\tau\choose 2}{n\choose 2}}
\prod_{1\le j<k\le n}
(z_jz_k)^\tau(z_j/z_k)_\tau(z_k/z_j)_\tau,
\label{eq:(z-z)tau}
\end{eqnarray}
the function $\Phi(z)\Delta(z)$ is symmetric. Therefore we obtain (\ref{eq:I(zeta)=I()}). $\square$
\begin{remark} As pointed out in Lemma \ref{lem:I-vanishing}, 
the right-hand side of (\ref{eq:I(zeta)=I()}) makes sense only when $\tau$ is a positive integer.  
However, as a function the left-hand side of (\ref{eq:I(zeta)=I()}) is defined continuously 
whether $\tau$ is a positive integer or not. 
Thus, as our basic strategy we first obtain several results for $I(x)$ 
under the condition $\tau\not\in\mathbb{Z}_+$. 
Then, the corresponding results in the cases $\tau\in\mathbb{Z}_+$ follow using
analytic continuation. 
Furthermore, if necessary they will be rewritten
using the relation (\ref{eq:I(zeta)=I()}), 
as we will see later. 
\end{remark}
%

We now state the formula corresponding to Proposition \ref{prop:aomoto}, extending $\tau$ from a positive integer 
to a complex number.

\begin{prop}[Aomoto] 
\label{thm:aomoto}
Under the condition {\rm (\ref{eq:convergence2-I})}, the truncated Jackson integral $I(\zeta)$ is expressed as  
\begin{equation}
\label{eq:I(zeta)2}
I(\zeta)=
(1-q)^n\prod_{j=1}^n (a_1t^{j-1})^{\alpha+2(n-j)\tau}
\frac{(q)_\8 (t)_\8
(q^{\alpha}a_1b_1t^{n+j-2})_\8
}
{(t^j)_\8(q^{\alpha}t^{j-1})_\8(a_1b_1t^{j-1})_\8}.
\end{equation}
\end{prop}
\noindent
{\bf Proof.} See the proof written after Lemma \ref{lem:I(a+N;zeta)} in the next subsection. 
As an alternative proof, see also Corollary \ref{cor:04} below. 
$\square$
\begin{remark} If we substitute $a_1$ and $b_1$ as $a_1\to 1$ and $b_1\to q^\beta$, respectively, 
then (\ref{eq:I(zeta)2}) coincides with (\ref{eq:I(zeta)}). 
\end{remark}

Moreover, we have the multiple bilateral sum version of the above formula as follows:
\begin{prop}[Aomoto] 
\label{thm:I(x)}
Suppose $\tau\not\in \mathbb{Z}_+$. 
Under the condition {\rm (\ref{eq:convergence1-I})}, $I(x)$ is expressed as a ratio of theta functions:
\begin{equation}
\label{eq:I(x)}
I(x)=c_0\, \prod_{i=1}^n x_i^{\alpha+2(n-i)\tau}
\frac{\theta(q^{\alpha}b_1t^{n-1}x_i)
}{\theta(b_1 x_i)}\prod_{1\le j<k\le n}\frac{\theta(x_k/x_j)}{\theta(t x_k/x_j)},
\end{equation}
where $c_0$ is a constant independent of $x\in (\mathbb{C}^*)^n$, which is explicitly written as 
\begin{equation}
\label{eq:c0}
c_0=(1-q)^n
\prod_{j=1}^{n}
\frac{(q)_\8(qt^{-j})_\8
(qa_1^{-1}b_1^{-1}t^{-(j-1)})_\8
}
{(qt^{-1})_\8(q^{\alpha}t^{j-1})_\8(q^{1-\alpha}a_1^{-1}b_1^{-1}t^{-(n+j-2)})_\8
}.
\end{equation}
\end{prop}
\begin{remark}
If $n=1$, then (\ref{eq:I(x)}) is equivalent to Ramanujan's $_1\psi_1$ summation theorem. 
This is another multi-dimensional bilateral extension of Ramanujan's $_1\psi_1$ summation theorem, 
which is different from the Milne--Gustafson summation theorem \cite{Gu87,Mi86}. 
(See also \cite{IF12} for the explanation of the Milne--Gustafson summation theorem along the context of this paper.) 
Another class of extension relates to the theory of Macdonald polynomials;
see for example \cite{Kan96,MS02,W05} as cited in \cite{W13}. 
\end{remark}
\noindent
{\bf Proof.}
Taking account of the poles of $\Phi(z)$, we have the expression 
\begin{equation}
\label{eq:I(x)poles}
I(x)=f(x)\prod_{i=1}^n \frac{x_i^\alpha}{\theta(b_1 x_i)}
\prod_{1\le j<k\le n}
\frac{x_j^{2\tau-1}}{\theta(t x_k/x_j)},
\end{equation}
where $f(x)$ is some holomorphic function on $(\mathbb{C}^*)^n$.
Under the condition $\tau\not\in \mathbb{Z}_+$, 
from Lemma \ref{lem:I-vanishing}, $I(x)$ is divisible by  
$
\prod_{1\le i<j\le n}
x_i\theta(x_j/x_i)
$.
This indicates that 
$
f(x)=g(x)\prod_{1\le i<j\le n}
x_i\theta(x_j/x_i),
$
where $g(x)$ is a holomorphic function on $(\mathbb{C}^*)^n$. 
Taking account of the $q$-periodicity of both sides of (\ref{eq:I(x)poles}), we have
$$
T_{x_i}g(x)=-\frac{g(x)}{q^{\alpha}b_1t^{n-1}x_i}
\quad\mbox{for}\quad i=1,2,\ldots, n,
$$
where $T_{x_i}$ means the shift operator of $x_i\to qx_i$, i.e.,
$T_{x_i}g(\ldots,x_i,\ldots)=g(\ldots,qx_i,\ldots)$. 
Then $g(x)$ is uniquely determined as
$
g(x)=c_0\prod_{i=1}^n \theta(q^{\alpha}b_1t^{n-1}x_i),
$
where $c_0$ is a constant independent of $x$. Therefore we obtain the expression (\ref{eq:I(x)}). 
Comparing (\ref{eq:I(zeta)2}) with (\ref{eq:I(x)}) of $x=\zeta$, the explicit form of $c_0$ is obtained as (\ref{eq:c0}). 
$\square$\\

Combining Lemma \ref{lem:I(zeta)=I()} and Proposition \ref{thm:I(x)},  
we obtain a multiple bilateral summation formula extending (\ref{eq:AHKE}) in Proposition \ref{prop:AHKE}. 
\begin{cor}
\label{cor:I(x,...,x)}
Suppose $\tau\in \mathbb{Z}_+$. For an arbitrary $x\in \mathbb{C}^*$
\begin{equation}
\label{eq:I(x,...,x)}
\frac{I(x,x,\ldots,x)}{n!}=c_1\, \prod_{i=1}^n (xt^{i-1})^{\alpha+2(n-i)\tau}
\frac{\theta(q^{\alpha}b_1t^{n+i-2}x)
}{\theta(b_1t^{i-1} x)}
,
\end{equation}
where $c_1$ is given by
\begin{equation}
\label{eq:c1}
c_1
=(1-q)^n
\prod_{j=1}^{n}
\frac{(q)_\8(t)_\8
(qa_1^{-1}b_1^{-1}t^{-(j-1)})_\8
}
{(t^{j})_\8(q^{\alpha}t^{j-1})_\8(q^{1-\alpha}a_1^{-1}b_1^{-1}t^{-(n+j-2)})_\8
}.
\end{equation}
\end{cor}
\noindent
{\bf Proof.} 
First we temporarily assume $\tau\not\in \mathbb{Z}_+$. 
From Proposition \ref{thm:I(x)}, we immediately have the following for the point $(x,xt,\ldots,xt^{n-1})$:
\begin{eqnarray*}
I(x,xt,\ldots,xt^{n-1})&=&c_0\, \prod_{i=1}^n (xt^{i-1})^{\alpha+2(n-i)\tau}
\frac{\theta(q^{\alpha}b_1t^{n+i-2}x)\theta(t)
}{\theta(b_1t^{i-1} x)\theta(t^{i})}\\
&=&c_1\, \prod_{i=1}^n (xt^{i-1})^{\alpha+2(n-i)\tau}
\frac{\theta(q^{\alpha}b_1t^{n+i-2}x)
}{\theta(b_1t^{i-1} x)},
\end{eqnarray*}
where $c_1$ is written as (\ref{eq:c1}). Then, by analytic continuation, the above formula is valid for $\tau\in \mathbb{Z}_+$. 
Using Lemma \ref{lem:I(zeta)=I()}, we obtain (\ref{eq:I(x,...,x)}). $\square$
\begin{remark} (\ref{eq:I(x,...,x)}) of Corollary \ref{cor:I(x,...,x)} is also expressed as 
$$
\frac{I(x,x,\ldots,x)}{n!}=
\frac{I(a_1,a_1,\ldots,a_1)}{n!}
\prod_{i=1}^n \Big(\frac{x}{a_1}\Big)^{\!\!\alpha +(n-1)\tau}
\frac{\theta(q^{\alpha}xb_1t^{n+i-2})\theta(a_1b_1t^{i-1})}{\theta(q^{\alpha}a_1b_1t^{n+i-2})\theta(xb_1t^{i-1})},
$$
which is the connection between $I(x,x,\ldots,x)$ and $I(a_1,a_1,\ldots,a_1)$. 
\end{remark}

As a special case of Corollary \ref{cor:I(x,...,x)}, we immediately have the formula (\ref{eq:AHKE}) in Proposition \ref{prop:AHKE}.
\begin{cor}[Askey, Habsieger, Kadell, Evans]
Suppose $\tau\in \mathbb{Z}_+$. Then 
$$\frac{I(a_1,a_1,\ldots,a_1)}{n!}
=(1-q)^n\prod_{j=1}^n (a_1t^{j-1})^{\alpha+2(n-j)\tau}
\frac{(q)_\8 (t)_\8
(q^{\alpha}a_1b_1t^{n+j-2})_\8
}
{(t^j)_\8(q^{\alpha}t^{j-1})_\8(a_1b_1t^{j-1})_\8}. 
$$
\end{cor}
\begin{remark} 
If we substitute $a_1$ and $b_1$  
as $a_1\to 1$ and $b_1\to q^\beta$, respectively, 
then the above formula coincides with (\ref{eq:AHKE}). 
\end{remark}

\subsection{$q$-difference equation with respect to $\alpha$}
In this subsection we derive the $q$-difference equation with respect to $\alpha$ satisfied by $I(x)$. 
We use $I(\alpha;x)$ instead of $I(x)$ to highlight the $\alpha$ dependence. 
The following lemma is known as Aomoto's method \cite{Ao87,AAR99,Fo10}. 
\begin{lem}[Aomoto]
\label{lem:aomoto}
Let $e_i(z), i=0,1,\ldots,n$, be the elementary symmetric polynomials, i.e., 
$$
e_r(z):=\sum_{1\le i_1< i_2<\cdots< i_r\le n}z_{i_1}z_{i_2}\cdots z_{i_r}
 \quad\mbox{for}\quad r=1,2,\ldots, n,
$$
and $e_0(z):=1$. Then 
\begin{equation}
\label{eq:Se=Se-Aomoto}
\int_0^{\mbox{\small $x$}\8} 
e_{i}(z)\Phi(z)\Delta(z)\varpi_q
=\frac{a_1t^{i-1}(1-t^{n-i+1})(1-q^{\alpha}t^{n-i})}
{(1-t^{i})(1-q^{\alpha}a_1b_1t^{2n-i-1})}
\int_0^{\mbox{\small $x$}\8} 
e_{i-1}(z)\Phi(z)\Delta(z)\varpi_q.
\end{equation}
\end{lem}
\noindent
{\bf Proof.} See \cite[Chapter 4, Exercises 4.6 q.2]{Fo10}. $\square$\\

Since we have 
$$
I(\alpha+1;x)=
\int_0^{\mbox{\small $x$}\8}z_1z_2\cdots z_n
\Phi(z)\Delta(z)\varpi_q,
$$
by definition, using (\ref{eq:Se=Se-Aomoto}) repeatedly we immediately obtain the $q$-difference equation with respect to $\alpha$ as
\begin{cor}
\label{cor:rec}
The recurrence relation for $I(\alpha;x)$ is given by 
\begin{equation}
\label{eq:01rec1}
I(\alpha+1;x)=I(\alpha;x)
\prod_{i=1}^n\frac{a_1t^{i-1}(1-q^{\alpha}t^{i-1})}{1-q^{\alpha}a_1b_1t^{n+i-2}}.
\end{equation}
\end{cor}

By definition 
$I(\zeta)$ is a sum over $\Lambda\cong{\mathbb N}^n$, while $I(x)$ is generally the sum over the lattice ${\mathbb Z}^n$.
It has the advantage of simplifying the computation of the $\alpha\to +\8$ asymptotic behavior,
as will be seen below.
(The lattice  $\{(x_1q^{\nu_1},\ldots,x_nq^{\nu_n})\in ({\mathbb C}^*)^n\,;\,(\nu_1,\ldots,\nu_n)\in {\mathbb Z}^n\}$ is 
called the {\it $q$-cycle} \cite{Ao91} of $I(x)$, while the set 
$\{(a_1q^{\nu_1},a_1tq^{\nu_2},\ldots,a_1t^{n-1}q^{\nu_n})\in ({\mathbb C}^*)^n\,;\,(\nu_1,\ldots,\nu_n)\in \Lambda\}$ as the support of $I(\zeta)$ is called 
the {\it $\alpha$-stable cycle} in \cite{Ao94,AK94-2}.)
\begin{lem}
\label{lem:I(a+N;zeta)}
The asymptotic behavior of $I(\alpha+N;\zeta)$ as $N\to +\8$ is given by
\begin{equation}
\label{eq:01asym1}
I(\alpha+N;\zeta)\sim (1-q)^n
\prod_{i=1}^n(a_1t^{i-1})^{\alpha+2(n-i)\tau+N}\frac{(q)_\8(t)_\8}{(a_1b_1t^{i-1})_\8(t^{i})_\8}
\quad(N\to +\8).
\end{equation}
\end{lem}
\noindent
{\bf Proof.} Since $\Phi(\zeta q^{\nu})=\Phi(a_1q^{\nu_1},a_1tq^{\nu_2},\ldots,a_1t^{n-1}q^{\nu_n})=0$ if $\nu\not\in \Lambda$, by definition $I(\zeta)$ is written as 
\begin{eqnarray*}
I(\alpha+N;\zeta)&=&(1-q)^n
\sum_{0\le\nu_1\le \nu_2\le \cdots\le \nu_n}
\prod_{i=1}^n(a_1t^{i-1}q^{\nu_i})^{\alpha+2(n-i)\tau+N}
\frac{(t^{i-1}q^{1+\nu_i})_\8}{(a_1b_1t^{i-1}q^{\nu_i})_\8}\\
&&\hskip 70pt\times
\prod_{1\le j<k\le n}\frac{(t^{k-j-1}q^{1+\nu_k-\nu_j})_\8}{(t^{k-j+1}q^{\nu_k-\nu_j})_\8}
(1-t^{k-j}q^{\nu_k-\nu_j}),
\end{eqnarray*}
so that the leading term of the asymptotic behavior of $I(\alpha+N;\zeta)$ as $N\to +\8$ is given by the 
term corresponding to $(\nu_1,\ldots,\nu_n)=(0,\ldots,0)$ in the above sum, which is (\ref{eq:01asym1}).
$\square$ \\

\noindent
{\bf Proof of Proposition \ref{thm:aomoto}.} By repeated use of the recurrence relation (\ref{eq:01rec1}), we have 
$$I(\alpha;x)
=I(\alpha+N;x)
\prod_{i=1}^n\frac{(q^{\alpha}a_1b_1t^{n+i-2})_N}
{(a_1t^{i-1})^N(q^{\alpha}t^{i-1})_N}
.$$
If we put $x=\zeta$ and take $N\to +\8$, we obtain  
\begin{equation}
\label{eq:I(alpha;zeta)}
I(\alpha;\zeta)
=
\lim_{N\to \8}\frac{I(\alpha+N;\zeta)}
{\prod_{i=1}^n(a_1t^{i-1})^N}
\times
\prod_{i=1}^n
\frac{(q^{\alpha}a_1b_1t^{n+i-2})_\8}
{(q^{\alpha}t^{i-1})_\8},
\end{equation}
which coincides with the right-hand side of (\ref{eq:I(zeta)2}) if we use (\ref{eq:01asym1}).
This means that the truncated Jackson integral $I(\zeta)$ is the special solution of the $q$-difference equation 
(\ref{eq:01rec1}), fixed by the asymptotic behavior (\ref{eq:01asym1}) as $\alpha\to +\infty$. 
\subsection{Regularization and connection formula}
Let ${\cal I}(x)$ and $h(x)$ be the functions defined by 
\begin{equation}
\label{eq:01cal I(x)1}
{\cal I}(x)=\frac{I(x)}{h(x)} \quad\mbox{where}\quad h(x):=\prod_{i=1}^n \frac{x_i^\alpha}{\theta(b_1 x_i)}
\prod_{1\le j<k\le n}x_j^{2\tau}\frac{\theta(x_k/x_j)}{\theta(t x_k/x_j)}.
\end{equation}
We call ${\cal I}(x)$ the {\it regularized Jackson integral} of $I(x)$. 
Since the trivial poles and zeros of $I(x)$ are canceled out by multiplying together $1/h(x)$ and $I(x)$, 
we have the following. 
\begin{lem}
\label{lem:h-and-s}
The regularization ${\cal I}(x)$ is holomorphic on $({\mathbb C}^*)^n$ and symmetric. 
\end{lem}
\noindent
{\bf Proof.} From the expression (\ref{eq:Phi0}) of $\Phi(z)$ as integrand of (\ref{eq:I(xi)}), 
the function $I(x)$ has the poles lying only in the set 
$\{x=(x_1,x_2,\ldots,x_n)\in ({\mathbb C}^*)^n\,;\,\prod_{i=1}^n\theta(b_1 x_i)
\prod_{1\le i<j\le n}\theta(t x_j/x_i)=0\}$.
Moreover, from Lemma \ref{lem:I-vanishing}, $I(x)$ is divisible by $x_j\theta(x_i/x_j)$. We therefore obtain 
$$
I(x)={\cal I}(x)h(x),
$$
where ${\cal I}(x)$ is some holomorphic function on $({\mathbb C}^*)^n$. 
Since $h(x)$ also satisfies $\sigma h(x)=(\sgn\,\sigma)U_\sigma(x)h(x)$ as (\ref{eq:Usigma-I}), 
${\cal I}(x)$ is symmetric. $\square$\\

From Proposition \ref{thm:I(x)} the regularization ${\cal I}(x)$ is written as 
\begin{equation}
\label{eq:01cal I(x)2}
{\cal I}(x)=
c_0\prod_{i=1}^n \theta(q^{\alpha}b_1t^{n-1}x_i).
\end{equation}
\begin{lem}[connection formula]
For an arbitrary $x,y\in ({\mathbb C}^*)^n$, the {\it connection formula} between ${\cal I}(x)$ and ${\cal I}(y)$ 
is written as 
\begin{equation}
\label{eq:01cal I(x)I(y)}
{\cal I}(x)={\cal I}(y)\prod_{i=1}^n
\frac{\theta(q^{\alpha}b_1t^{n-1}x_i)}{\theta(q^{\alpha}b_1t^{n-1}y_i)}.
\end{equation}
In particular, if we set $y=\zeta\in ({\mathbb C}^*)^n$, then 
\begin{equation}
\label{eq:01cal I(x)I(a)}
{\cal I}(x)={\cal I}(\zeta)\prod_{i=1}^n
\frac{\theta(q^{\alpha}b_1t^{n-1}x_i)}{\theta(q^{\alpha}a_1b_1t^{n+i-2})}.
\end{equation}
\end{lem}
\noindent
{\bf Proof.}
From (\ref{eq:01cal I(x)2}), we immediately have (\ref{eq:01cal I(x)I(y)}), 
and recalling (\ref{3.4'}) gives (\ref{eq:01cal I(x)I(a)}). $\square$
\begin{remark} 
In the equation (\ref{eq:01cal I(x)I(a)}), if we switch the symbols from ${\cal I}(x)$ to $I(x)$ we obtain 
\begin{equation}
\label{eq:01I(x)I(a)}
{I}(x)=I(\zeta)\frac{h(x)}{h(\zeta)}
\prod_{i=1}^n
\frac{\theta(q^{\alpha}b_1t^{n-1}x_i)}{\theta(q^{\alpha}a_1b_1t^{n+i-2})}, 
\end{equation}
which is also the connection formula between a solution $I(x)$ of the $q$-difference equation (\ref{eq:01rec1}) 
and the special solution $I(\zeta)$ fixed by its 
asymptotic behavior (\ref{eq:01asym1}) as $\alpha\to +\8$. In addition, its connection coefficient is written as 
a ratio of theta functions (i.e., that of $q$-gamma functions), and is of course invariant under the shift 
$\alpha\to \alpha+1$.  
From the evaluation (\ref{eq:I(zeta)2}) of $I(\zeta)$, 
the connection formula (\ref{eq:01I(x)I(a)}) is another expression for 
the product formula (\ref{eq:I(x)}) in Proposition \ref{thm:I(x)}.
\end{remark}
If we set 
\begin{equation}
\label{eq:01beta}
\beta:=1-\alpha_1-\beta_1-2(n-1)\tau-\alpha, 
\end{equation}
where $\alpha_1$ and $\beta_1$ are given by $a_1=q^{\alpha_1}, b_1=q^{\beta_1}$, 
after rearrangement, 
the formula (\ref{eq:I(x)}) is also expressed as the following Macdonald-type sum, 
whose value is given by an $x$-independent constant 
\cite{Ma03,vD97,Ito06-1}.
\begin{prop} Under the condition $a_1b_1t^{2n-2}q^{\alpha+\beta}=q$, 
\begin{eqnarray}
\label{eq:01cal I(x)4}
&&\!\!\!\!\!\!
\int_0^{\mbox{\small $x$}\8}\prod_{i=1}^n
\frac{(qa_1^{-1}z_i)_\infty(qb_1^{-1}z_i^{-1})_\infty}
{(q^{\alpha}b_1t^{n-1}z_i)_\8(q^{\beta}a_1t^{n-1}z_i^{-1})_\8}
\prod_{1\le j<k\le n}
\frac{(qt^{-1}z_j/z_k)_\8(qt^{-1}z_k/z_j)_\8}{(qz_j/z_k)_\8(qz_k/z_j)_\8}
\,\frac{d_qz_1}{z_1}\wedge\cdots\wedge\frac{d_qz_n}{z_n}
\nonumber\\
&&
\qquad
=(1-q)^n
\prod_{j=1}^{n}
\frac{(q)_\8(qt^{-j})_\8
(qa_1^{-1}b_1^{-1}t^{-(j-1)})_\8
}
{(qt^{-1})_\8(q^{\alpha}t^{j-1})_\8(q^{\beta}t^{j-1})_\8
}.
\end{eqnarray}
\end{prop} 
\noindent
{\bf Proof.} Since $h(x)\prod_{i=1}^n\theta(q^{\alpha}b_1t^{n-1}x_i)$ is invariant under the $q$-shift $x_i\to qx_i$, 
from (\ref{eq:01cal I(x)2}), we have 
\begin{equation*}
\int_0^{\mbox{\small $x$}\8}\frac{\Phi(z)\Delta(z)}
{h(z)\prod_{i=1}^n\theta(q^{\alpha}b_1t^{n-1}z_i)}\,\frac{d_qz_1}{z_1}\wedge\cdots\wedge\frac{d_qz_n}{z_n}
=c_0,
\end{equation*}
so that 
\begin{equation*}
\int_0^{\mbox{\small $x$}\8}\prod_{i=1}^n
\frac{(qa_1^{-1}z_i)_\infty(qb_1^{-1}z_i^{-1})_\infty}
{\theta(q^{\alpha}b_1t^{n-1}z_i)}
\prod_{1\le j<k\le n}
\frac{(qt^{-1}z_j/z_k)_\8(qt^{-1}z_k/z_j)_\8}{(qz_j/z_k)_\8(qz_k/z_j)_\8}\,\frac{d_qz_1}{z_1}\wedge\cdots\wedge\frac{d_qz_n}{z_n}
=c_0,
\end{equation*}
which is rewritten as (\ref{eq:01cal I(x)4}) using (\ref{eq:c0}) under the condition (\ref{eq:01beta}). $\square$\\

As a corollary, it is confirmed that the following identity for a contour integral is equivalent to  
the formula (\ref{eq:01cal I(x)4}) of the special case $x=\zeta$.    

\begin{cor}
Let $\mathbb{T}^n$ be the the direct product of the unit circle, i.e.,  
$
\mathbb{T}^n:=\{(z_1,\ldots,z_n)\in\mathbb{C}^n \,;\,|z_i|=1\}
$. 
Suppose that $|a_1|<1$, $|b_1|<1$, $|t|<1$ and $a_1b_1t^{2n-2}q^{\alpha+\beta}=q$. 
Then
\begin{eqnarray}
\label{eq:contour}
&&\!\!\!\!\!\!
\Big(\frac{1}{2\pi\sqrt{-1}}\Big)^{\!\!n}\frac{1}{n!}
\int_{\mathbb{T}^n}
\prod_{i=1}^n \frac{(q^\alpha t^{n-1}a_1z_i^{-1})_\8(q^\beta t^{n-1}b_1 z_i)_\8}{(a_1z_i^{-1})_\8(b_1 z_i)_\8}
\prod_{1\le j<k\le n}\frac{(z_j/z_k)_\8(z_k/z_j)_\8}{(t z_j/z_k)_\8(t z_k/z_j)_\8}
\frac{dz_1}{z_1}\cdots\frac{dz_n}{z_n}\nonumber\\
&&\quad=\prod_{i=1}^n \frac{(t)_\infty(q^{1-\beta}t^{-(i-1)})_\infty(q^{1-\alpha}t^{-(i-1)})_\infty}
{(q)_\infty(t^i)_\infty(a_1b_1t^{i-1})_\infty}.
\end{eqnarray}
\end{cor}
\noindent
{\bf Proof.} This follows by a residue calculation using (\ref{eq:01cal I(x)4}) of the case $x=\zeta$. $\square$

%

%
\subsection{Dual expression of the Jackson integral $I(x)$}
For an arbitrary $x=(x_1,x_2,\ldots,x_n)\in ({\mathbb C}^*)^n$ 
we specify $x^{-1}$ as 
\begin{equation}
\label{eq:01x-1}
x^{-1}:=(x_1^{-1},x_2^{-1},\ldots,x_n^{-1})\in ({\mathbb C}^*)^n.
\end{equation}
For the point 
$\bar\zeta=(b_1,b_1t,\ldots,b_1t^{n-1})\in ({\mathbb C}^*)^{n}$, 
if we set $y=\bar\zeta^{-1}$ in the connection formula (\ref{eq:01cal I(x)I(y)}), 
then we obtain the expression 
\begin{equation}
\label{eq:01cal I(x)I(b)}
{\cal I}(x)={\cal I}(\bar\zeta^{-1})\prod_{i=1}^n
\frac{\theta(q^{\alpha}b_1t^{n-1}x_i)}{\theta(q^{\alpha}t^{i-1})}.
\end{equation}
Since $x=\bar\zeta^{-1}=(b_1^{-1},b_1^{-1}t^{-1},\ldots,b_1^{-1}t^{-(n-1)})$ is a pole of the function $I(x)$ by definition,  
$I(\bar\zeta^{-1})$ no longer makes sense. 
However, the regularization ${\cal I}(\bar\zeta^{-1})$ appearing on the right-hand side of (\ref{eq:01cal I(x)I(b)}) still has 
meaning as a special value of a holomorphic function. 
We will show a way to realize the regularization ${\cal I}(\bar\zeta^{-1})$ as a computable object by another Jackson integral.  
For this purpose, let $\bar\Phi(z)$ be the function specified by 
\begin{equation}
\label{eq:barPhi}
\bar\Phi(z):=
\prod_{i=1}^nz_i^{1-\alpha_1-\beta_1-2(n-1)\tau-\alpha}\frac{(qb_1^{-1}z_i)_\8}{(a_1z_i)_\8}
\prod_{1\le j<k\le n}
z_j^{2\tau-1}\frac{(qt^{-1}z_k/z_j)_\8}{(t z_k/z_j)_\8},
\end{equation}
where $\alpha_1$ and $\beta_1$ are given by $a_1=q^{\alpha_1}, b_1=q^{\beta_1}$. 
For $x=(x_1,x_2,\ldots,x_n)\in ({\mathbb C}^*)^n$, we define the sum $\bar I(x)$ by  
\begin{equation}
\label{eq:01bar I(x)}
\bar I(x):=\int_0^{\mbox{\small $x$}\8}
\bar\Phi(z)\Delta(z)\varpi_q,
\end{equation}
which converges absolutely under the condition (\ref{eq:convergence1-I}).
We call $\bar I(x)$ the {\it dual Jackson integral of $I(x)$}, 
and call $\bar I(\bar\zeta)$ its {\it truncation}. 
Note that the sum $I(x)$ transforms to its dual $\bar I(x)$ if we interchange the parameters as 
\begin{equation}
\label{eq:transform}
\alpha\leftrightarrow\beta
\quad\mbox{and}\quad a_1\leftrightarrow b_1,
\end{equation}
where $\beta$ is specified by (\ref{eq:01beta}). 
We also define the {\it regularization} $\bar {\cal I}(x)$ of $\bar I(x)$ as
\begin{equation}
\label{eq:01barh(x)}
\bar{\cal I}(x)=\frac{\bar I(x)}{\bar h(x)} \quad\mbox{where}\quad 
\bar h(x):=\prod_{i=1}^n \frac{x_i^{1-\alpha_1-\beta_1-2(n-1)\tau-\alpha}}{\theta(a_1x_i)}
\prod_{1\le j<k\le n}x_j^{2\tau}
\frac{\theta(x_k/x_j)}{\theta(t x_k/x_j)}.
\end{equation}
In the same manner as Lemma \ref{lem:h-and-s}, we can confirm that the function $\bar {\cal I}(x)$ is also holomorphic and symmetric.

\begin{lem}[reflective equation]
\label{lem:01ref}
The connection between $I(x)$ and $\bar I(x)$ is 
\begin{equation}
\label{eq:01ref1}
I(x)=\frac{h(x)}{\bar h(x^{-1})}\bar I(x^{-1}),
\end{equation}
where $x^{-1}$ is specified as in {\rm (\ref{eq:01x-1})} and 
$$
\frac{h(x)}{\bar h(x^{-1})}=
(-1)^{n\choose 2}
\prod_{i=1}^n x_i^{1-\alpha_1-\beta_1}\frac{\theta(qa_1^{-1}x_i)}{\theta(b_1 x_i)}
\prod_{1\le j<k\le n}\Big(\frac{x_k}{x_j}\Big)^{\!\!1-2\tau}\frac{\theta(qt^{-1}x_k/x_j)}{\theta(t x_k/x_j)}.
$$
In other words, the relation between ${\cal I}(x)$ and $\bar {\cal I}(x)$ is 
\begin{equation}
\label{eq:01ref2}
{\cal I}(x)=\bar{\cal I}(x^{-1}).
\end{equation}
\end{lem}
\noindent
{\bf Proof.} 
From the definitions (\ref{eq:01cal I(x)1}) and (\ref{eq:01barh(x)}) the ratio $h(x)/\bar h(x^{-1})$ 
is written as in (\ref{eq:01ref1}). 
Since $
\Delta(z)=(-1)^{{n\choose 2}}(z_1z_2\cdots z_n)^{n-1}\Delta(z^{-1})
$, 
from (\ref{eq:Phi0}), (\ref{eq:Delta}), (\ref{eq:barPhi}), 
we have 
\begin{equation}
\label{eq:01PD=hhPD}
\Phi(z)\Delta(z)=
\frac{h(z)}{\bar h(z^{-1})}
\bar\Phi(z^{-1})\Delta(z^{-1}).
\end{equation}
Also since $h(z)/\bar h(z^{-1})$ is invariant under the shift $z_i\to qz_i$,   
by the definitions (\ref{eq:I(xi)}) and (\ref{eq:01bar I(x)}) of $I(z)$ and $\bar I(z)$, 
the connection (\ref{eq:01ref1}) between $I(z)$ and its dual $\bar I(z)$ is derived from (\ref{eq:01PD=hhPD}). $\square$\\

We use $\bar I(\alpha;x)$ instead of $\bar I(x)$ to see the $\alpha$ dependence. 
From (\ref{eq:01ref1}), the recurrence relation for $\bar I(\alpha;x)$ is completely the same as (\ref{eq:01rec1}) of $I(\alpha;x)$. 
\begin{lem}
\label{lem:rec2}
The function $\bar I(\alpha;x)$ also satisfies the recurrence relation {\rm (\ref{eq:01rec1})} of $I(\alpha;x)$, and is rewritten as  
\begin{equation}
\label{eq:01rec2}
{\bar I}(\alpha;x)={\bar I}(\alpha-1;x)
\prod_{i=1}^n\frac{1-q^{1-\alpha}t^{-(i-1)}}
{b_1t^{i-1}(1-q^{1-\alpha}a_1^{-1}b_1^{-1}t^{-(n+i-2)})}.
\end{equation}
\end{lem}

We saw above that although $I(\bar\zeta^{-1})$ no longer makes sense,  
its regularization ${\cal I}(\bar\zeta^{-1})$ still has meaning as a special value of a holomorphic function, 
and ${\cal I}(\bar\zeta^{-1})$ is evaluated by the dual integral $\bar{\cal I}(\bar\zeta)$ 
via the reflective equation (\ref{eq:01ref2}). Moreover, by definition the regularization $\bar{\cal I}(\bar\zeta)$ itself is calculated 
using $\bar I(\bar\zeta)$, which is then normally defined as a truncated Jackson integral.  
Though we already know the value of ${\cal I}(\bar\zeta^{-1})$ 
through the connection formula (\ref{eq:01cal I(x)I(a)}) with $x=\bar\zeta^{-1}$, 
the point is that 
we can calculate ${\cal I}(\bar\zeta^{-1})$ directly from $\bar I(\bar\zeta)$. 
This requires the leading term of its asymptotic behavior as $\alpha\to -\8$
is simply computed as follows. 

\begin{cor}
\label{cor:bar I(a-N;b)}
The asymptotic behavior of $\bar I(\alpha-N;\bar\zeta)$ as $N\to +\8$ is written as
\begin{eqnarray}
\label{eq:01asym2}
\bar I(\alpha-N;\bar\zeta)
\sim (1-q)^n
\prod_{i=1}^n(b_1t^{i-1})^{
1-\alpha_1-\beta_1-2(i-1)\tau-\alpha+N}\frac{(q)_\8(t)_\8}{(a_1b_1t^{i-1})_\8(t^{i})_\8}
\quad(N\to +\8).
\end{eqnarray}
Moreover, by repeated use of {\rm (\ref{eq:01rec2})}, the truncated Jackson integral $\bar I(\bar\zeta)$ is written as  
\begin{equation}
\label{eq:01I(b)}
\bar I(\bar\zeta)=
(1-q)^n\prod_{i=1}^n 
\frac{(b_1t^{i-1})^{1-\alpha_1-\beta_1-2(i-1)\tau-\alpha}(q)_\8 (t)_\8
(q^{1-\alpha}t^{-(i-1)})_\8
}
{(t^i)_\8(q^{1-\alpha}a_1^{-1}b_1^{-1}t^{-(n+i-2)})_\8(a_1b_1t^{i-1})_\8}.
\end{equation}
\end{cor}
\noindent
{\bf Proof.}
Using Lemma \ref{lem:rec2}, the arguments are completely parallel to 
Lemma \ref{lem:I(a+N;zeta)} and (\ref{eq:I(alpha;zeta)}). 
Actually, using the rule (\ref{eq:transform}), if we substitute $a_1$, $b_1$ and $\alpha$ in $\Phi(z)$ of (\ref{eq:Phi0})
by $b_1$, $a_1$ and $\beta$, respectively, then $\Phi(z)$ transforms to $\bar\Phi(z)$ in (\ref{eq:barPhi}), 
so that we obtain the same result as Lemma \ref{lem:I(a+N;zeta)} and Proposition \ref{thm:aomoto} with these substitutions. $\square$\\

From (\ref{eq:01cal I(x)I(y)}) and  (\ref{eq:01ref2}), for $x, y\in ({\mathbb C}^*)^{n}$ 
we have the connection formula between ${\cal I}(x)$ and $\bar{\cal I}(y)$ as 
$$
{\cal I}(x)=\bar{\cal I}(y)\prod_{i=1}^n
\frac{\theta(q^{\alpha}b_1t^{n-1}x_i)}{\theta(q^{\alpha}b_1t^{n-1}y_i^{-1})}.
$$
In particular, if $y=\bar\zeta$, then we have 
$$
{\cal I}(x)=\bar{\cal I}(\bar\zeta)
\prod_{i=1}^n
\frac{\theta(q^{\alpha}b_1t^{n-1}x_i)}{\theta(q^{\alpha}t^{i-1})}. 
$$
If we switch the symbols from ${\cal I}(x)$ and $\bar{\cal I}(\bar\zeta)$ to $I(x)$ and $\bar I(\bar\zeta)$, respectively, 
then we obtain 
\begin{equation}
\label{eq:01I(x)I(b)}
I(x)={\bar I}(\bar\zeta)
\frac{h(x)}{{\bar h}(\bar\zeta)}
\prod_{i=1}^n
\frac{\theta(q^{\alpha}b_1t^{n-1}x_i)}{\theta(q^{\alpha}t^{i-1})}
.
\end{equation}
We once again obtain the connection formula between a solution $I(x)$ of 
the $q$-difference equation (\ref{eq:01rec1}) and the special solution $\bar I(\bar\zeta)$ fixed by its 
asymptotic behavior (\ref{eq:01asym2}) as $\alpha\to -\8$, 
as a counterpart of the formula (\ref{eq:01I(x)I(a)}) of the case $\alpha\to +\8$.

The connection formula (\ref{eq:01I(x)I(b)}) with (\ref{eq:01I(b)}) is also another expression 
for the product formula (\ref{eq:I(x)}) in Proposition \ref{thm:I(x)}, 
like the formula (\ref{eq:01I(x)I(a)}). 

\begin{remark}
The truncated Jackson integrals $I(\zeta)$ and ${\bar I}(\bar\zeta)$ both satisfy the $q$-difference equation (\ref{eq:01rec1}) with respect to $\alpha$. $I(\zeta)$ is the special solution fixed by the asymptotic behavior (\ref{eq:01asym1}) as $\alpha\to +\infty$. On the other hand, ${\bar I}(\bar\zeta)$ is the solution fixed by the asymptotic behavior (\ref{eq:01asym2}) as $\alpha\to -\infty$. The connection formula (\ref{eq:01I(x)I(b)}) shows  
\begin{equation}
I(\zeta)={\bar I}(\bar\zeta)
\prod_{i=1}^n
\frac{(a_1t^{i-1})^{\alpha+2(n-i)\tau}\theta(q^{\alpha}a_1b_1t^{n+i-2})}
{(b_1t^{i-1})^{1-\alpha_1-\beta_1-2(i-1)\tau-\alpha}\theta(q^{\alpha}t^{i-1})},
\end{equation}
which connects $I(\zeta)$ and ${\bar I}(\bar\zeta)$ by the $q$-periodic function of the right-hand side.
This formula is explained like the formula $\Gamma(\alpha)\Gamma(1-\alpha)=\pi/\sin \pi \alpha$, which indicates that 
$\Gamma(\alpha)$ and $1/\Gamma(1-\alpha)$ are solutions of the difference equation
$f(\alpha+1)=\alpha f(\alpha)$ and they are fixed by the specific asymptotic behaviors 
(i.e., Stirling's formula) as $\alpha\to +\infty$ and $-\infty$, respectively, and these solutions are connected by 
the periodic function $\pi/\sin \pi \alpha$.
\end{remark}

\noindent
\begin{remark}
As we have seen above, we used the integrand $\bar\Phi(z)$ instead of $\Phi(z)$, 
which coincides with $\bar\Phi(z)$
up to the $q$-periodic factor $h(z)/h(z^{-1})$, and used the set 
$\{(b_1q^{\nu_1},b_1tq^{\nu_2},\ldots,b_1t^{n-1}q^{\nu_n})$\\$\in ({\mathbb C}^*)^n\,;\,(\nu_1,\ldots,\nu_n)\in \Lambda\}$
as the ``$(-\alpha)$-stable cycle" for the dual integral $\bar I(x)$  
when we construct a special solution $\bar I(\bar\zeta)$ expressed by 
(Jackson) integral representation for the $q$-difference equation (\ref{eq:01rec1}) as $\alpha\to -\8$ . 
In the classical setting, this process is usually done by taking an imaginary cycle without changing the integrand $\Phi(z)$
under the ordinary integral representation. In the $q$-analog setting  Aomoto and Aomoto--Kato 
used the integral representation without changing the integrand $\Phi(z)$, but instead, 
they adopted the residue sum on the set 
$\{(b_1^{-1}q^{-\nu_1},b_1^{-1}t^{-1}q^{-\nu_2},\ldots,b_1^{-1}t^{-(n-1)}q^{-\nu_n})\in ({\mathbb C}^*)^n\,;\,(\nu_1,\ldots,\nu_n)\in \Lambda\}$
of poles of $I(x)$. They call this cycle the {\it $\alpha$-unstable cycle} \cite{Ao94,AK94-2} of $I(x)$ for the parameter $\alpha$.
To carry out this process is called {\it regularization} in their original paper \cite{Ao90}.
We hope our slight changes of terminology does not bring confusion to the reader. 
\end{remark}

\section{Further extension of the $q$-Selberg integral}
\label{section:02}
In this section we will explain a bilateral extension of Askey--Evans's $q$-Selberg integral.
\subsection{Jackson integral of Selberg type}
\label{subsection:02.1}
Let $a_1,a_2,b_1,b_2$ and $t$ be complex numbers satisfying 
\begin{equation}
\label{eq:convergence1}
|qt^{2i-2}|<1\quad\mbox{and}\quad q<|a_1a_2b_1b_2t^{2i-2}|\quad\mbox{for}\quad i=1,2,\ldots,n.
\end{equation}
Let $\Phi(z)$ be specified by 
\begin{equation}
\label{eq:Phi1}
\Phi(z):=\prod_{i=1}^n z_i\frac{(qa_1^{-1}z_i)_\8}{(b_1z_i)_\8}\frac{(qa_2^{-1}z_i)_\8}{(b_2z_i)_\8}
\prod_{1\le j<k\le n}z_j^{2\tau-1}\frac{(qt^{-1}z_k/z_j)_\8}{(tz_k/z_j)_\8},
\end{equation}
where $\tau$ is given by $t=q^\tau$, and let $\Delta(z)$ be the difference product specified by (\ref{eq:Delta}). 
For $x=(x_1,x_2,\ldots,x_n)\in (\mathbb{C}^*)^n$, we define the sum $J(x)$ by
\begin{equation}
\label{eq:J(x)}
J(x)=J(x_1,x_2,\ldots,x_n):=
\int_0^{\mbox{\small $x$}\8}
\Phi(z)\Delta(z)
\varpi_q,
\end{equation}
which converges absolutely under the condition (\ref{eq:convergence1}).
We call $J(x)$ the {\it Jackson integral of Selberg type}. 
For arbitrary $x_1,x_2\in \mathbb{C}^*$, we set the points $\zeta_i(x_1,x_2)\in (\mathbb{C}^*)^n$ by
\begin{equation}
\label{eq:zeta(x1x2)}
\zeta_i(x_1,x_2):=(\hskip 1pt\underbrace{x_1,x_1t,\ldots,x_1t^{i-1}}_{i},
\underbrace{x_2,x_2t,\ldots,x_2t^{n-i-1}}_{n-i}\hskip 1pt)\in (\mathbb{C}^*)^n
\quad\mbox{for}\quad
i=0,1,\ldots,n.
\end{equation}
When we set $x_1=a_1$ and $x_2=a_2$ on (\ref{eq:zeta(x1x2)}), for the special points $\zeta_i(a_1,a_2)$, $i=0,1,\ldots,n$, 
by definition $J(\zeta_i(a_1,a_2))$ is defined as the sum (\ref{eq:J(x)}) over the fan region 
\begin{equation}
\label{eq:Lambda-i}
\Lambda_i=\{(\nu_1,\nu_2,\ldots,\nu_n)\in \mathbb{Z}^n\,;\,
0\le \nu_1\le \nu_2\le \cdots\le \nu_i\ \mbox{and}\  0\le \nu_{i+1}\le \nu_{i+2}\le \ldots\le \nu_n\}.
\end{equation}
We call $J(\zeta_i(a_1,a_2))$
the {\em truncated Jackson integral of Selberg type}. \\

The main theorem of this paper is the following:
\begin{thm} 
\label{thm:main1}
\begin{equation}
\label{eq:main1}
\sum_{i=0}^{n}J(\zeta_i(x_1,x_2))
\prod_{j=1}^{i}\frac{(x_1t^{j-1})^{-1-2(n-j)\tau}}{\theta(x_2x_1^{-1}t^{-j+1})}
\prod_{k=1}^{n-i}\frac{(x_2t^{k-1})^{-1-2(n-i-k)\tau}}{\theta(x_1x_2^{-1}t^{i-k+1})}
=C_0\prod_{k=1}^n
\frac
{\theta(x_1x_2b_1b_2t^{n+k-2})}
{\prod_{i=1}^2\prod_{j=1}^2\theta(x_ib_jt^{k-1})},
\end{equation}
where $C_0$ is a constant independent of $x_1$ and $x_2$. The constant $C_0$ is explicitly expressed as  
\begin{equation}
\label{eq:C0}
C_0=(1-q)^n\prod_{k=1}^n\frac{(q)_\8(t)_\8\prod_{i=1}^2\prod_{j=1}^2(qa_i^{-1}b_j^{-1}t^{-(k-1)})_\8}
{(t^{k})_\8(qa_1^{-1}a_2^{-1}b_1^{-1}b_2^{-1}t^{-(n+k-2)})_\8}.
\end{equation}
\end{thm}
\begin{remark} In particular, if $\tau\in\mathbb{Z}_+$,  then the formula simplifies
\begin{equation}
\label{eq:main1remark}
\sum_{i=0}^{n}(-1)^iJ(\zeta_i(x_1,x_2))=C_0
(-1)^{\tau{n\choose 2}}q^{-{\tau\choose 2}{n\choose 2}}
\prod_{j=1}^n
\frac
{(x_1x_2t^{j-1})^{(n-j)\tau}x_2\theta(x_1/x_2)\theta(x_1x_2b_1b_2t^{n+j-2})}
{\theta(x_1b_1t^{j-1})\theta(x_1b_2t^{j-1})\theta(x_2b_1t^{j-1})\theta(x_2b_2t^{j-1})},
\end{equation}
which is equivalent to Corollary \ref{cor:main1+} as will be explained later.
\end{remark}

We will give the proof of Theorem \ref{thm:main1} in the next subsection. 
(See the proof of Theorem \ref{thm:main2}, which is equivalent to Theorem \ref{thm:main1}.)
In this subsection we will explain the relation between this main theorem and other known results 
which are deduced from this theorem as corollaries. As a special case $x_1=a_1$, $x_2=a_2$ of Theorem \ref{thm:main1}, 
we have the formula for the truncated Jackson integrals $J(\zeta_i(a_1,a_2))$.
\begin{cor} [Tarasov--Varchenko, Stokman]
\label{cor:main1}
\begin{eqnarray}
\label{eq:TV,S}
&&\sum_{i=0}^{n}J(\zeta_i(a_1,a_2))
\prod_{j=1}^{i}\frac{(a_1t^{j-1})^{-1-2(n-j)\tau}}{\theta(a_2a_1^{-1}t^{-j+1})}
\prod_{k=1}^{n-i}\frac{(a_2t^{k-1})^{-1-2(n-i-k)\tau}}{\theta(a_1a_2^{-1}t^{i-k+1})}
\nonumber\\
&&\quad
=(1-q)^n\prod_{k=1}^n
\frac
{(q)_\8(t)_\8(a_1a_2b_1b_2t^{n+k-2})_\8}
{(t^{k})_\8\prod_{i=1}^2\prod_{j=1}^2(a_ib_jt^{k-1})_\8}.
\end{eqnarray}
\end{cor}
\begin{remark} 
This formula was given by Tarasov--Varchenko \cite[Theorem (E.10)]{TV97} and Stokman \cite[Corollary 7.6]{St00} independently. 
The proof of \cite{TV97} is by computing residues for an A type generalization of Askey--Roy's $q$-beta integral, 
while the proof of \cite{St00} is by computing residues for Gustafson's $q$-Selberg contour integral \cite{Gu90,Gu94} 
and an appropriate limiting procedure. 
This formula extends  Askey--Evans's formula (\ref{eq:Evans0}) from $\tau$ a positive integer to $\tau$ an arbitrary complex number.
(Compare with Corollary \ref{cor:main1++}.) 
\end{remark} 

As we saw in the section of the Jackson integral of A-type before, and in keeping with the above remark, 
it is also very important for the Jackson integral of Selberg type 
to distinguish between the situations of the cases 
whether the parameter $\tau$ is a positive integer or not.  

\begin{lem}
\label{lem:J=0}
Suppose $\tau\not\in\mathbb{Z}_+$. 
If $x_i=x_j$ for some $i$ and $j$ $(1\le i<j\le n)$, then $J(x_1,x_2,\ldots,x_n)=0$. 
\end{lem}
\noindent
{\bf Proof.} In the same way as Lemma \ref{lem:I-vanishing}. $\square$\\

On the other hand, under the condition $\tau$ being a positive integer 
we generally have $J(x_1,\ldots,x_n)\ne 0$ even if $x_i=x_j$ $(1\le i<j\le n)$. 
In particular, we have the following:

\begin{lem}
\label{lem:J(zeta)=J()}
Suppose $\tau\in\mathbb{Z}_+$. 
For arbitrary $x_1,x_2\in \mathbb{C}^*$, $J(\zeta_i(x_1,x_2))$ is expressed as 
\begin{equation}
\label{eq:J(zeta)=J()}
J(\zeta_i(x_1,x_2))
=\frac{1}{i!(n-i)!}J(\hskip 1pt
\underbrace{x_1,x_1,\ldots,x_1}_{i},
\underbrace{x_2,x_2,\ldots,x_2}_{n-i}
\hskip 1pt). 
\end{equation}
\end{lem}
\noindent
{\bf Proof.} In the same way as Lemma \ref{lem:I(zeta)=I()}. $\square$
\begin{remark} As pointed out in Lemma \ref{lem:J=0}, 
the right-hand side of (\ref{eq:J(zeta)=J()}) makes sense only when $\tau$ is a positive integer.  
However, as a function of $\tau$ the left-hand side of (\ref{eq:J(zeta)=J()}) is defined continuously 
whether $\tau$ is a positive integer or not. 
Thus, as our basic strategy we first obtain several results for $J(x)$ 
under the condition $\tau\not\in\mathbb{Z}_+$. 
Then, using analytic continuation, 
the results for $J(x)$ can automatically be regarded as those for $\tau\in\mathbb{Z}_+$. 
And if necessary, we will rewrite them appropriately using the relation (\ref{eq:J(zeta)=J()}). 
\end{remark}

Recalling the binomial theorem
$$
(x_2-x_1)^n=\sum_{i=0}^n(-1)^i{n\choose i}x_2^{n-i}x_1^i, 
$$
under the condition $\tau\in\mathbb{Z}_+$  and using (\ref{eq:00jac3}) and Lemma \ref{lem:J(zeta)=J()}, 
we can deform the following iterated Jackson integral
as
\begin{eqnarray}
\label{eq:SSSddd}
&&\frac{1}{n!}\int_{x_1\infty}^{x_2\infty}\hskip -10pt\cdots\int_{x_1\infty}^{x_2\infty}\hskip -2pt
\Phi(z)\Delta(z)
\,\frac{d_qz_1}{z_1}\cdots\frac{d_qz_n}{z_n}
\nonumber\\[1pt]
&&=\frac{1}{n!}\sum_{i=0}^n(-1)^i{n\choose i}
\underbrace{\int_{0}^{x_2\infty}\hskip -10pt\cdots\int_{0}^{x_2\infty}\hskip -16pt}_{n-i}
\hskip 7pt 
\underbrace{\int_{0}^{x_1\infty}\hskip -10pt\cdots\int_{0}^{x_1\infty}\hskip -16pt}_i
\hskip 14pt 
\Phi(z)\Delta(z)
\,\frac{d_qz_1}{z_1}\cdots\frac{d_qz_n}{z_n}
\nonumber\\[-7pt]
&&=\sum_{i=0}^n(-1)^i\frac{1}{i!(n-i)!}
J(\hskip 1pt 
\underbrace{x_1,x_1,\ldots,x_1}_{i},
\underbrace{x_2,x_2,\ldots,x_2}_{n-i}
\hskip 1pt )
\nonumber\\
%
&&=\sum_{i=0}^n(-1)^i
J(\zeta_i(x_1,x_2)),
\end{eqnarray}
because the integrand $\Phi(z)\Delta(z)$ is symmetric if $\tau\in\mathbb{Z}_+$. 
 \\
 
If $\tau\in \mathbb{Z}_+$, using (\ref{eq:00quasi-period2}), 
the coefficient factor of $J(\zeta_i(x_1,x_2))$ appearing in the left-hand side of (\ref{eq:main1}) is simplified as 
\begin{eqnarray*}
&&
(-1)^{i}
\prod_{j=1}^{i}(x_1t^{j-1})^{1+2(n-j)\tau}\theta(x_2x_1^{-1}t^{-j+1})
\prod_{k=1}^{n-i}(x_2t^{k-1})^{1+2(n-i-k)\tau}\theta(x_1x_2^{-1}t^{i-k+1})\\
&&\quad\quad
=(-1)^{\tau{n\choose 2}}q^{-{\tau\choose 2}{n\choose 2}+\tau^2{n\choose 3}}
(x_1x_2)^{\tau{n\choose 2}}
\big(x_2\theta(x_1/x_2)\big)^n,
\end{eqnarray*}
which is independent of the choice of indices $i=0,1,\ldots,n$. 
From this, if $\tau\in\mathbb{Z}_+$, 
the formula (\ref{eq:main1}) in the main theorem then shrinks to the form (\ref{eq:main1remark}).
From (\ref{eq:SSSddd}), (\ref{eq:main1remark}) is rewritten to 
\begin{eqnarray*}
&&\frac{1}{n!}\int_{x_1\infty}^{x_2\infty}\hskip -10pt\cdots\int_{x_1\infty}^{x_2\infty}\hskip -2pt
\Phi(z)\Delta(z)
\,\frac{d_qz_1}{z_1}\cdots\frac{d_qz_n}{z_n}
\\
&&\quad\quad
=(-1)^{\tau {n \choose 2}}q^{-{\tau \choose 2}{n \choose 2}}C_0\prod_{j=1}^n
\frac
{(x_1x_2t^{j-1})^{(n-j)\tau}x_2\theta(x_1/x_2)\theta(x_1x_2b_1b_2t^{n+j-2})}
{\theta(x_1b_1t^{j-1})\theta(x_1b_2t^{j-1})\theta(x_2b_1t^{j-1})\theta(x_2b_2t^{j-1})}.
\end{eqnarray*}
Using (\ref{eq:(z-z)tau}), we therefore obtain
\begin{cor}
\label{cor:main1+} Suppose $\tau\in\mathbb{Z}_+$. Then 
\begin{eqnarray}
\label{eq:Evans3}
&&\frac{1}{n!}\int_{x_1\infty}^{x_2\infty}\!\!\!\cdots\int_{x_1\infty}^{x_2\infty}
\prod_{i=1}^n
z_i^{(n-1)\tau}
\frac{(qa_1^{-1}z_i)_\8}{(b_1z_i)_\8}
\frac{(qa_2^{-1}z_i)_\8}{(b_2z_i)_\8}
\prod_{1\le j<k\le n}
(z_j/z_k)_\tau(z_k/z_j)_\tau
\,d_qz_1\cdots d_qz_n
\nonumber\\
&&\hskip 10mm=C_0
\prod_{j=1}^n
\frac{(x_1x_2q^{(j-1)\tau})^{(n-j)\tau}x_2\theta(x_1/x_2)\,\theta(x_1x_2b_1b_2q^{(n+j-2)\tau})}
{\theta(x_1b_1q^{(j-1)\tau})\theta(x_1b_2q^{(j-1)\tau})\theta(x_2b_1q^{(j-1)\tau})\theta(x_2b_2q^{(j-1)\tau})},
\end{eqnarray}
where $C_0$ is the constant given by {\rm (\ref{eq:C0})}.\end{cor}

In particular, putting $x_1=a_1$ and $x_2=a_2$ on the above equation, we obtain
\begin{cor}
[Askey, Evans]
\label{cor:main1++} 
Suppose $\tau\in\mathbb{Z}_+$. Then 
\begin{eqnarray}
\label{eq:Evans4}
&&\frac{1}{n!}\int_{a_1}^{a_2}\cdots\int_{a_1}^{a_2}
\prod_{i=1}^n
z_i^{(n-1)\tau}
\frac{(qa_1^{-1}z_i)_\8}{(b_1z_i)_\8}
\frac{(qa_2^{-1}z_i)_\8}{(b_2z_i)_\8}
\prod_{1\le j<k\le n}
(z_j/z_k)_\tau(z_k/z_j)_\tau
\,d_qz_1\cdots d_qz_n
\nonumber\\
&&\hskip 10pt=(1-q)^n\prod_{j=1}^n
\frac{(q)_\infty(t)_\infty(a_1a_2b_1b_2t^{n+j-2})_\infty\, (a_1a_2t^{j-1})^{(n-j)\tau}a_2\theta(a_1/a_2)}
{(t^{j})_\infty(a_1b_1t^{j-1})_\infty(a_1b_2t^{j-1})_\infty(a_2b_1t^{j-1})_\infty(a_2b_2t^{j-1})_\infty}.
\end{eqnarray}
\end{cor}
\begin{remark}
If we substitute $a_1, a_2, b_1,b_2$  
as $a_1\to x_1$, $a_2\to x_2$, $b_1\to q^\alpha/x_1 $, $b_2\to q^\beta/x_2$, respectively, 
then (\ref{eq:Evans4}) is rewritten as
\begin{eqnarray*}
&&\frac{1}{n!}\int_{x_1}^{x_2}\cdots\int_{x_1}^{x_2}
\prod_{i=1}^n
z_i^{(n-1)\tau}
\frac{(qz_i/x_1)_\8}{(q^{\alpha}z_i/x_1)_\8}\frac{(qz_i/x_2)_\8}{(q^{\beta}z_i/x_2)_\8}
\prod_{1\le j<k\le n}
(z_j/z_k)_\tau(z_k/z_j)_\tau
\,d_qz_1\cdots d_qz_n
\nonumber\\
&&\hskip 17pt=
\prod_{j=1}^n\frac{(1-q)(q)_\8(q^{\alpha+\beta+(n+j-2)\tau})_\8(q^\tau)_\8}
{(q^{\alpha+(j-1)\tau})_\8(q^{\beta+(j-1)\tau})_\8
(q^{j\tau})_\8}
\frac{(x_1x_2q^{(j-1)\tau})^{(n-j)\tau}x_2\theta(x_1/x_2)}{(x_2q^{\alpha+(j-1)\tau}/x_1)_\8(x_1q^{\beta+(j-1)\tau}/x_2)_\8},
\end{eqnarray*}
which exactly coincides with the formula (\ref{eq:Evans0}) established by Askey and Evans.
\end{remark}

\subsection{Regularization and the proof of the main theorem}
Let ${\cal J}(x)$ be the function defined by 
\begin{equation}
\label{eq:calJ}
{\cal J}(x):=\frac{J(x)}{h(x)}
\quad\mbox{where}\quad
h(x):=\prod_{i=1}^n\frac{x_i}{\theta(b_1x_i)\theta(b_2x_i)}
\prod_{1\le j<k\le n}x_j^{2\tau}\frac{\theta(x_k/x_j)}{\theta(tx_k/x_j)}.
\end{equation}
\begin{lem}
The function ${\cal J}(x)$ is holomorphic on $(\mathbb{C}^*)^n$ and symmetric. 
\end{lem}
\noindent
{\bf Proof.} In the same way as Lemma \ref{lem:h-and-s}. $\square$\\

We call ${\cal J}(x)$ the {\em regularization of $J(x)$} or the {\em regularized Jackson integral of $J(x)$}.
For the point $\zeta_i(x_1,x_2)$ defined by (\ref{eq:zeta(x1x2)}),  $h(\zeta_i(x_1,x_2))$ is evaluated as 
\begin{equation}
\label{eq:h(zeta(x1x2))}
h(\zeta_i(x_1,x_2))=
\prod_{j=1}^i\frac{(x_1t^{j-1})^{1+2(n-j)\tau}}{\theta(b_1x_1t^{j-1})\theta(b_2x_1t^{j-1})}
\frac{\theta(x_2x_1^{-1}t^{-(j-1)})}{\theta(x_2x_1^{-1}t^{n-i-j+1})}\frac{\theta(t)}{\theta(t^j)}
\prod_{k=1}^{n-i}\frac{(x_2t^{k-1})^{1+2(n-i-k)\tau}}{\theta(b_2x_1t^{k-1})\theta(b_2x_2t^{k-1})}
\frac{\theta(t)}{\theta(t^k)},
\end{equation}
which is used below. \\

Using the definition (\ref{eq:calJ}) of the regularization ${\cal J}(x)$, 
it is directly confirmed that 
the equation (\ref{eq:main1}) in the main theorem (Theorem \ref{thm:main1}) is 
rewritten in the following form:
\begin{thm}
\label{thm:main2}
Suppose $\tau\not\in\mathbb{Z}_+$. 
Let $H_i(x_1,x_2)$ be the function defined as 
$$
H_i(x_1,x_2):=
\prod_{j=1}^{i}\frac{\theta(x_2b_1t^{n-j})\theta(x_2b_2t^{n-j})}{\theta(t^{j})\theta(x_2x_1^{-1}t^{n-i-j+1})}
\prod_{k=1}^{n-i}\frac{\theta(x_1b_1t^{n-k})\theta(x_1b_2t^{n-k})}{\theta(t^{k})\theta(x_1x_2^{-1}t^{i-k+1})}
\prod_{l=1}^n\frac{\theta(t^l)}{\theta(x_1x_2b_1b_2t^{n+l-2})}.
$$
Then 
\begin{equation}
\label{eq:main2-1}
\sum_{i=0}^{n}{\cal J}(\zeta_i(x_1,x_2))H_i(x_1,x_2)=C_1,
\end{equation}
where $C_1$ is a constant independent of $x_1$ and $x_2$.
\end{thm}
\begin{remark} The constant $C_1$ is explicitly given by 
\begin{equation}
\label{eq:C1}
C_1=(1-q)^n\prod_{k=1}^n\frac{(q)_\8(qt^{-k})_\8\prod_{i=1}^2\prod_{j=1}^2(qa_i^{-1}b_j^{-1}t^{-(k-1)})_\8}
{(qt^{-1})_\8(qa_1^{-1}a_2^{-1}b_1^{-1}b_2^{-1}t^{-(n+k-2)})_\8},
\end{equation}
which will be confirmed later. 
\end{remark}

From Theorem \ref{thm:main2} we immediately have an expression for the constant $C_1$.
\begin{cor}
\label{cor:C1=J}
The constant $C_1$ in {\rm (\ref{eq:main2-1})} is expressed by the special values of ${\cal J}(x)$ as 
\begin{equation}
\label{eq:C1=J}
C_1={\cal J}(\zeta_i(b_1^{-1}t^{-(i-1)},b_2^{-1}t^{-(n-i-1)}))
\quad(i=0,1,\ldots,n). 
\end{equation}
\end{cor}
\noindent
{\bf Proof.} 
Since $H_i(x_1,x_2)$ has the property 
$$
H_i(b_1^{-1}t^{-(j-1)},b_2^{-1}t^{-(n-j-1)})=\delta_{ij}\quad(i,j=0,1,\ldots,n),
$$
where $\delta_{ij}$ is Kronecker's delta, using (\ref{eq:main2-1}), 
we obtain the constant $C_1$ as (\ref{eq:C1=J}). $\square$
\begin{remark} Since $x=\zeta_i(b_1^{-1}t^{-(i-1)},b_2^{-1}t^{-(n-i-1)})$ is a pole of the function $J(x)$ by definition,  
the value $J(\zeta_i(b_1^{-1}t^{-(i-1)},b_2^{-1}t^{-(n-i-1)}))$ no longer makes sense. 
However, the regularization ${\cal J}(\zeta_i(b_1^{-1}t^{-(i-1)},b_2^{-1}t^{-(n-i-1)}))$ still has 
meaning as a special value of a holomorphic function. In the next subsection 
we will show a way to realize the regularization ${\cal J}(\zeta_i(b_1^{-1}t^{-(i-1)},b_2^{-1}t^{-(n-i-1)}))$ as a computable object by another Jackson integral.  
And eventually it will lead us to the explicit evaluation of the constant $C_1$ as (\ref{eq:C1}), 
which will be confirmed later as Lemmas \ref{lem:C1=barJ} and \ref{lem:barJn}. 
\end{remark}

The rest of this subsection is devoted to the proof of Theorem \ref{thm:main2}. 
If we set 
\begin{equation}
\label{eq:F(x1,x2)}
F(x_1,x_2)=\sum_{i=0}^{n}F_i(x_1,x_2), 
\end{equation}
where
$$
F_i(x_1,x_2):={\cal J}(\zeta_i(x_1,x_2))
\prod_{j=1}^{i}\frac{\theta(x_2b_1t^{n-j})\theta(x_2b_2t^{n-j})}{\theta(t^{j})\theta(x_2x_1^{-1}t^{n-i-j+1})}
\prod_{k=1}^{n-i}\frac{\theta(x_1b_1t^{n-k})\theta(x_1b_2t^{n-k})}{\theta(t^{k})\theta(x_1x_2^{-1}t^{i-k+1})}
\prod_{l=1}^n\theta(t^l), 
$$
then 
the equation (\ref{eq:main2-1}) is equivalent to  
\begin{equation}
\label{eq:main2-2}
F(x_1,x_2)=C_1\prod_{j=1}^n\theta(x_1x_2b_1b_2t^{n+j-2}). 
\end{equation}
In order to prove Theorem \ref{thm:main2} we will show (\ref{eq:main2-2}) instead of (\ref{eq:main2-1}). For this purpose, we will prove two lemmas first. 
\begin{lem}
\label{lem:dimH=n}
Suppose that $\tau\not\in\mathbb{Z}_+$. 
Let $H$ be the set of the holomorphic functions on $(\mathbb{C}^*)^2$ satisfying 
the $q$-difference equation
\begin{equation}
\label{eq:quasi-period}
f(qx_1,x_2)=f(x_1,qx_2)=(-x_1x_2b_1b_2t^{3(n-1)/2})^{-n}f(x_1,x_2).
\end{equation}
The dimension of $H$ as a linear space is equal to $n$, i.e., $\dim_{\mathbb{C}} H=n$.
Moreover, the set $\{\Theta_i(x_1,x_2)\,;\, i=1,2,\ldots,n\}$ is a basis of $H$, where 
$\Theta_i(x_1,x_2)$ is defined by 
\begin{equation}
\label{eq:Theta-i(x1,x2)}
\Theta_i(x_1,x_2):=\theta(x_1x_2b_1b_2t^{n-i})
\prod_{1\le j\le n\atop j\ne i}\theta(x_1x_2b_1b_2t^{2n-j}),\quad i=1,2,\ldots,n. 
\end{equation}
\end{lem}
\begin{remark} $\Theta_1(x_1,x_2)$ coincides with the function appearing in the right-hand side of (\ref{eq:main2-2}).  
\end{remark}
\noindent
{\bf Proof.} For an arbitrary function $f(x_1,x_2)\in H$, 
since $f(x_1,x_2)$ is holomorphic function on $(\mathbb{C}^*)^2$, 
$f(x_1,x_2)$ is expanded as 
$
f(x_1,x_2)=\sum_{i,j=-\infty}^\infty c_{ij}x_1^{i}x_2^{j}.
$
From $f(qx_1,x_2)=f(x_1,qx_2)$, we have $c_{ij}q^i=c_{ij}q^j$. 
This indicates that $c_{ij}=0$ if $i\ne j$.
Denoting $c_{ii}$ by $c_i$, $f(x_1,x_2)$ is written as 
$
f(x_1,x_2)=\sum_{i=-\infty}^\infty c_{i}(x_1x_2)^{i}.
$
From $f(qx_1,x_2)=(-x_1x_2b_1b_2t^{3(n-1)/2})^{-n}f(x_1,x_2)$, we have 
$c_iq^i=c_{i+n}(-b_1b_2t^{3(n-1)/2})^{-n}$. 
This indicates that $f(x_1,x_2)$ is determined by $c_0,c_1,\ldots,c_{n-1}$, 
which means that $H$ is spanned by its $n$ elements.

Since it is obvious that $\Theta_i(x_1,x_2)\in H$ from the explicit expression (\ref{eq:Theta-i(x1,x2)}), 
it suffices for our purpose 
to show the linearly independence of $\{\Theta_i(x_1,x_2)\}$. Assume that 
$\sum_{i=1}^nc_i\Theta_i(x_1,x_2)=0$. By definition 
\begin{equation}
\label{eq:Theta-Kronecker}
\Theta_i(b_1^{-1}t^{-(n-j+1)},b_2^{-1}t^{-(n-1)})
=\delta_{ij}\,\theta(t^{-n})\prod_{1\le k\le n\atop k\ne j}\theta(t^{j-k})\quad (i,j=1,2,\ldots,n),
\end{equation}
where $\delta_{ij}$ is Kronecker's delta. 
Therefore we have $0=\sum_{i=1}^nc_i\Theta_i(b_1^{-1}t^{-(n-j+1)},b_2^{-1}t^{-(n-1)})=
c_j\,\theta(t^{-n})\prod_{1\le k\le n\atop k\ne j}\theta(t^{j-k}),$
which indicates $c_1=c_2=\cdots=c_n=0$. $\square$
\begin{lem} 
\label{lem:FinH}
Suppose 
$\tau\not\in\mathbb{Z}_+$.
Then 
$F(x_1,x_2)\in H$. 
\end{lem}
\noindent
{\bf Proof.} 
Since it is easy to confirm that $F(x_1,x_2)$ satisfies the $q$-difference equations (\ref{eq:quasi-period}), 
it suffices to prove that $F(x_1,x_2)$ is holomorphic on $(\mathbb{C}^*)^2$.  
For this purpose, in the expression (\ref{eq:F(x1,x2)}) of $F(x_1,x_2)$, 
we will confirm that the residues at the apparent poles vanish. 
Since each $F_i(x_1,x_2)$ in (\ref{eq:F(x1,x2)}) has the common quasi-periodicity (\ref{eq:quasi-period}), 
it suffices to consider the residues at the poles of the cases
(1) $x_1x_2^{-1}t^{i-j+1}=1$ $(1\le j\le n-i \le n)$ or (2) $x_2x_1^{-1}t^{n-i-j+1}=1$ $(1\le j\le i\le n)$.
We will examine these poles carefully below.
The function $F_0(x_1,x_2)$ has $n$ poles of order 1 at $x_1=x_2t^{j-1}$ $(1\le j\le n)$.
The function $F_n(x_1,x_2)$ also has $n$ poles of order 1 at $x_2=x_1t^{j-1}$ $(1\le j\le n)$.
If $i\ne 0,n$, then 
the function $F_i(x_1,x_2)$ is supposed to have $n$ poles at $x_1 t^{i}=x_2t^{j-1}$ $(1\le j\le n-i)$ and 
$x_2t^{n-i}=x_1t^{j-1}$ $(1\le j\le i)$. But when $x_1 t^{i}=x_2t^{j-1}$, 
if $1< j\le n-i$, from (\ref{eq:h(zeta(x1x2))}) and Lemma \ref{lem:J=0}, we have ${\cal J}(\zeta_i(x_1,x_2))=0$, which is a factor of $F_i(x_1,x_2)$. 
In the same manner, when $x_2t^{n-i}=x_1t^{j-1}$, 
if $1< j\le i$, we also have ${\cal J}(\zeta_i(x_1,x_2))=0$. 
This indicates that $F_i(x_1,x_2)$, $i\ne 0,n$, actually has 
only $2$ poles of order 1, i.e., $x_1 t^{i}=x_2$ and $x_2t^{n-i}=x_1$. Therefore 
the function $F(x_1,x_2)$ may have $2n-1$ poles of order 1 at $x_1=x_2t^{n-i}$ $(1\le i\le n)$ and 
$x_2=x_1t^{i}$ $(0\le i\le n-1)$ in total. 
However the residues there all vanish 
as is confirmed by the following calculation, valid for $1\le i\le n$:
\begin{eqnarray*}
&&
\!\!\!\!\!\!\!\!\!\!\!\!
\Res_{x_1=x_2t^{n-i}}F(x_1,x_2)=\lim_{x_1\to x_2t^{n-i}}(x_1-x_2t^{n-i})F(x_1,x_2)\\
&&=\lim_{x_1\to x_2t^{n-i}}(x_1-x_2t^{n-i})\big(F_0(x_1,x_2)+F_i(x_1,x_2)\big)\\
&&=\lim_{x_1\to x_2t^{n-i}}\frac{(x_1-x_2t^{n-i})}{\theta(x_1x_2^{-1}t^{-(n-i)})}
{\cal J}(\zeta_0(x_1,x_2))
\frac{\prod_{j=1}^{n}\theta(x_1b_1t^{n-j})\theta(x_1b_2t^{n-j})}{
\prod_{1\le j\le n\atop j\ne n-i+1}\theta(x_1x_2^{-1}t^{-j+1})}\\
&&\quad+\lim_{x_1\to x_2t^{n-i}}\frac{(x_1-x_2t^{n-i})}{\theta(x_2x_1^{-1}t^{n-i})}
{\cal J}(\zeta_i(x_1,x_2))
\prod_{j=1}^{n-i}\frac{\theta(x_1b_1t^{n-j})\theta(x_1b_2t^{n-j})}{\theta(t^{j})\theta(x_1x_2^{-1}t^{i-j+1})}\\
&&\hskip 40mm\times
\frac{\prod_{j=1}^{i}\theta(x_2b_1t^{n-j})\theta(x_2b_2t^{n-j})}{\prod_{j=2}^{i}\theta(t^{j})\theta(x_2x_1^{-1}t^{n-i-j+1})}
\prod_{j=2}^n\theta(t^j)\\
&&=-\frac{x_2t^{n-i}}{(q)_\infty^2}
\bigg[{\cal J}(\zeta_0(x_2t^{n-i},x_2))
-{\cal J}(\zeta_i(x_2t^{n-i},x_2))\bigg]
\frac{\prod_{j=1}^{n}\theta(x_2b_1t^{2n-i-j})\theta(x_2b_2t^{2n-i-j})}{
\prod_{1\le j\le n\atop j\ne n-i+1}\theta(t^{n-i-j+1})}\\
&&=0
\end{eqnarray*}
because ${\cal J}(\zeta_0(x_2t^{n-i},x_2))={\cal J}(\zeta_i(x_2t^{n-i},x_2))={\cal J}(x_2,x_2t,\ldots,x_2t^{n-1})$. 
In the same way as above it is also confirmed that 
$$\Res_{x_2=x_1t^i}F(x_1,x_2)
=\lim_{x_2\to x_1t^i}(x_2-x_1t^i)\big(F_i(x_1,x_2)+F_n(x_1,x_2)\big)
=0$$ for $0\le i\le n-1$.
Therefore $F(x_1,x_2)$ is holomorphic on $(\mathbb{C}^*)^2$. 
$\square$
\\[10pt]
\noindent
{\bf Proof of Theorem \ref{thm:main2}.} We will prove (\ref{eq:main2-2}).
From Lemmas \ref{lem:dimH=n} and \ref{lem:FinH}, 
$F(x_1,x_2)$ is expressed as a linear combination of $\Theta_i(x_1,x_2)$, $i=1,\ldots,n$, 
i.e.,
$
F(x_1,x_2)=\sum_{i=1}^nC_i\Theta_i(x_1,x_2), 
$
where $C_i$ are some constants. 
By definition, it is easy to confirm that
$F(b_1^{-1}t^{-(n-j+1)},b_2^{-1}t^{-(n-1)})$ $=0$ for $j=2,3,\ldots,n$. 
From (\ref{eq:Theta-Kronecker}) we therefore obtain 
$$
0=F(b_1^{-1}t^{-(n-j+1)},b_2^{-1}t^{-(n-1)})=\sum_{i=1}^nC_i\Theta_i(b_1^{-1}t^{-(n-j+1)},b_2^{-1}t^{-(n-1)})
=C_j\,\theta(t^{-n})\prod_{1\le k\le n\atop k\ne j}\theta(t^{j-k})
$$
for $j=2,3,\ldots,n$.
This indicates that $C_2=C_3=\cdots=C_n=0$. Thus we obtain (\ref{eq:main2-2}). $\square$

\subsection{Dual expression of the Jackson integral $J(x)$}
Let $\bar\Phi(z)$ be specified by 
\begin{equation}
\label{eq:Phi2-1}
\bar\Phi(z):=
\prod_{i=1}^n
z_i^{1-\alpha_1-\alpha_2-\beta_1-\beta_2-2(n-1)\tau}
\frac{(qb_1^{-1}z_i)_\8}{(a_1z_i)_\8}\frac{(qb_2^{-1}z_i)_\8}{(a_2z_i)_\8}
\prod_{1\le j<k\le n}
z_j^{2\tau-1}\frac{(qt^{-1}z_k/z_j)_\8}{(tz_k/z_j)_\8}, 
\end{equation}
where $\alpha_i$ and $\beta_i$ are given by $a_i=q^{\alpha_i}, b_i=q^{\beta_i}$, 
and let $\Delta(z)$ be specified by (\ref{eq:Delta}). 
For $x=(x_1,x_2,\ldots,x_n)\in (\mathbb{C}^*)^n$, we define the sum $\bar J(x)$ by
\begin{equation}
\label{eq:barJ(x)}
\bar J(x):=
\int_0^{\mbox{\small $x$}\8}
\bar\Phi(z)\Delta(z)\varpi_q,
\end{equation}
which converges absolutely under the condition (\ref{eq:convergence1}).
We call $\bar J(x)$ the {\it dual Jackson integral of $J(x)$}. 
When we set $x_1=b_1$ and $x_2=b_2$ in (\ref{eq:zeta(x1x2)}), for the special points $\zeta_i(b_1,b_2)$, $i=0,1,\ldots,n$, 
$\bar J(\zeta_i(b_1,b_2))$ is defined as the sum over the fan region $\Lambda_i$ specified by (\ref{eq:Lambda-i}).
We call $\bar J(\zeta_i(b_1,b_2))$
the {\em truncation} of the dual Jackson integral $\bar J(x)$. \\

Let $\bar{\cal J}(x)$ and $\bar h(x)$ be the functions defined by 
\begin{equation}
\label{eq:barcalJ}
\bar{\cal J}(x):=\frac{\bar J(x)}{\bar h(x)}
\quad\mbox{where}\quad
\bar h(x)=\prod_{i=1}^n\frac{x_i^{1-\alpha_1-\alpha_2-\beta_1-\beta_2-2(n-1)\tau}}{\theta(a_1x_i)\theta(a_2x_i)}
\prod_{1\le j<k\le n}x_j^{2\tau}\frac{\theta(x_k/x_j)}{\theta(tx_k/x_j)}.
\end{equation}
Since the trivial poles and zeros of $\bar J(x)$ are canceled out by multiplying together $1/\bar h(x)$ and $\bar J(x)$,
$\bar{\cal J}(x)$ is holomorphic on $x\in (\mathbb{C}^*)^n$ and symmetric. 
We call $\bar{\cal J}(x)$ the {\it regularization of $\bar J(x)$}. 
\begin{lem}[reflective equation]
\label{lem:ref}
For $x\in (\mathbb{C}^*)^n$, the relation between ${\cal J}(x)$ and $\bar{\cal J}(x)$ is given by 
\begin{equation}
\label{eq:ref3}
\bar{\cal J}(x)={\cal J}(x^{-1}),
\end{equation}
where $x^{-1}$ is specified as in {\rm (\ref{eq:01x-1})}. In particular, for $x_1,x_2\in \mathbb{C}^*$ 
the following holds:
\begin{equation}
\label{eq:ref4}
\bar{\cal J}(\zeta_i(x_1,x_2))={\cal J}(\zeta_i(x_1^{-1}t^{-(i-1)},x_2^{-1}t^{-(n-i-1)}))
\quad(i=0,1,\ldots,n). 
\end{equation}
\end{lem}
\noindent
{\bf Proof.} By definition (\ref{eq:ref3}) is equivalent to 
the connection between $J(x)$ and $\bar J(x)$, i.e., 
\begin{equation}
\label{eq:ref1}
\bar J(x)=\frac{\bar h(x)}{h(x^{-1})}J(x^{-1}),
\end{equation}
which we should prove. From the definitions (\ref{eq:calJ}) and (\ref{eq:barcalJ}) the ratio $\bar h(x)/h(x^{-1})$ 
is written as
\begin{equation}
\label{eq:ref2}
\frac{\bar h(x)}{h(x^{-1})}=
(-1)^{n\choose 2}
\prod_{i=1}^n\bigg(\prod_{l=1}^2 x_i^{1-\alpha_l-\beta_l}
\frac{\theta(qb_l^{-1}x_i)}{\theta(a_lx_i)}\bigg)
\prod_{1\le j<k\le n}\Big(\frac{x_k}{x_j}\Big)^{\!\!1-2\tau}\frac{\theta(qt^{-1}x_k/x_j)}{\theta(tx_k/x_j)}.
\end{equation}
Since $
\Delta(z)=(-1)^{n\choose 2}(z_1z_2\cdots z_n)^{n-1}\Delta(z^{-1})
$, 
from (\ref{eq:Phi1}), (\ref{eq:Phi2-1}) and (\ref{eq:Delta}), 
we have 
\begin{equation}
\label{eq:PD=hhPD}
\bar\Phi(z)\Delta(z)=
\frac{\bar h(z)}{h(z^{-1})}
\Phi(z^{-1})\Delta(z^{-1}).
\end{equation}
Also since $\bar h(z)/h(z^{-1})$ is invariant under the shift $z_i\to qz_i$,   
by the definitions (\ref{eq:J(x)}) and (\ref{eq:barJ(x)}) of $J(x)$ and $\bar J(x)$, 
the connection (\ref{eq:ref1}) between $J(x)$ and its dual $\bar J(x)$ is derived from (\ref{eq:PD=hhPD}). 
Since ${\cal J}(x)$ and $\bar{\cal J}(x)$ are symmetric, (\ref{eq:ref4}) is immediately followed from (\ref{eq:ref3}). $\square$\\

We now state the $q$-difference equations for $\bar J(x)$ under the setting 
$x=\zeta_i(b_1,b_2)$, $i=0,1,\ldots,n$. 

\begin{prop} 
\label{prop:q-diff01}
Suppose $x=\zeta_i(b_1,b_2), i=0,1,\ldots,n.$ Then 
the recurrence relations for $\bar J(x)$ are given by 
\begin{eqnarray}
T_{a_j}\bar J(x)&=&(-a_j)^n\prod_{k=1}^n
\frac{(1-a_j^{-1}b_1^{-1}t^{-(k-1)})(1-a_j^{-1}b_2^{-1}t^{-(k-1)})}{1-a_1^{-1}a_2^{-1}b_1^{-1}b_2^{-1}t^{-(n+k-2)}}\bar J(x),
\label{eq:TaJ=cJ}\\
T_{b_j}\bar J(x)&=&(-b_j^{-1})^n\prod_{k=1}^n
\frac{(1-a_1b_jt^{k-1})(1-a_2b_jt^{k-1})}{1-a_1a_2b_1b_2t^{n+k-2}}\bar J(x),
\label{eq:TbJ=cJ}
\end{eqnarray}
for $j=1,2$,
where $T_{a_j}$ and $T_{b_j}$ are the $q$-shift operators of $a_j\to qa_j$ and $b_j\to qb_j$, respectively. 
In other words, the recurrence relations for the regularization $\bar{\cal J}(x)$ are given by 
\begin{eqnarray}
T_{a_j}\bar {\cal J}(x)&=&\prod_{k=1}^n
\frac{(1-a_j^{-1}b_1^{-1}t^{-(k-1)})(1-a_j^{-1}b_2^{-1}t^{-(k-1)})}{1-a_1^{-1}a_2^{-1}b_1^{-1}b_2^{-1}t^{-(n+k-2)}}\bar {\cal J}(x),
\label{eq:TaJ=cJ2}\\
T_{b_j}\bar {\cal J}(x)&=&\prod_{k=1}^n
\frac{(1-a_1^{-1}b_j^{-1}t^{-(k-1)})(1-a_2^{-1}b_j^{-1}t^{-(k-1)})}{1-a_1^{-1}a_2^{-1}b_1^{-1}b_2^{-1}t^{-(n+k-2)}}\bar {\cal J}(x),
\label{eq:TbJ=cJ2}
\end{eqnarray}
for $j=1,2.$
\end{prop}
\noindent
{\bf Proof.} The derivation of (\ref{eq:TaJ=cJ}) and (\ref{eq:TbJ=cJ}) will be done in the Appendix. 
(See Remark \ref{remark:Se=Se} after Lemma \ref{lem:Se=Se}.) 
Here we just mention that (\ref{eq:TaJ=cJ2}) and (\ref{eq:TbJ=cJ2}) 
are derived from (\ref{eq:TaJ=cJ}) and (\ref{eq:TbJ=cJ}), respectively. 
From the expression (\ref{eq:barcalJ}) of $\bar h(x)$,  
under the condition $x=\zeta_i(b_1,b_2)$, $i=0,1,\ldots,n$, 
the function $\bar h(x)$ satisfies 
\begin{equation}
\label{eq:rec}
T_{a_j}\bar h(x)=(-a_j)^n\bar h(x)\quad\mbox{and}\quad
T_{b_j}\bar h(x)=\Big(\frac{b_j}{b_1b_2}\Big)^{\!n}t^{-{n\choose 2}}\,\bar h(x) 
\quad (j=0,1,\ldots, n).
\end{equation}
Since $\bar{\cal J}(x)={\bar J}(x)/{\bar h}(x)$, from the above equations and (\ref{eq:TaJ=cJ}), (\ref{eq:TbJ=cJ}) 
we therefore obtain (\ref{eq:TaJ=cJ2}), (\ref{eq:TbJ=cJ2}). $\square$

\subsection{Evaluation of the truncated Jackson integral}
The main result of this subsection is the evaluation of the regularization of the truncated Jackson integral 
using the $q$-difference equations (\ref{eq:TaJ=cJ2}) and (\ref{eq:TbJ=cJ2}) in Proposition \ref{prop:q-diff01} and its asymptotic behavior for the special direction of parameters. 
\begin{thm} 
\label{thm:main3}
For $x=\zeta_i(b_1,b_2), i=0,1,\ldots,n,$ 
the regularization $\bar{\cal J}(x)$ is evaluated as
$$
\bar{\cal J}(\zeta_i(b_1,b_2))=
(1-q)^n\prod_{k=1}^n\frac{(q)_\8(qt^{-k})_\8\prod_{i=1}^2\prod_{j=1}^2(qa_i^{-1}b_j^{-1}t^{-(k-1)})_\8}
{(qt^{-1})_\8(qa_1^{-1}a_2^{-1}b_1^{-1}b_2^{-1}t^{-(n+k-2)})_\8}.
$$
\end{thm}
This theorem can be deduced from the specific case of $i=n$ (or $i=0$). The reason is explained as follows. 
Using the reflective equation (\ref{eq:ref4}) and Corollary \ref{cor:C1=J}, we immediately have 
\begin{lem}
\label{lem:C1=barJ}
The constant $C_1$ in {\rm (\ref{eq:main2-1})} is expressed by the special values of the regularized Jackson integral as
$$
C_1=
\bar{\cal J}(\zeta_i(b_1,b_2)) \quad(i=0,1,\ldots,n). 
$$
\end{lem}
This indicates that $\bar{\cal J}(\zeta_i(b_1,b_2))$ does not depend on the choice of indices $i=0,1,\ldots,n$. 
From this fact, for the proof of Theorem \ref{thm:main3} it suffices to show that of the case $i=n$ only, i.e.,
\begin{lem}
\label{lem:barJn}
For $x=\bar\zeta=(b_1,b_1t,\ldots,b_1t^{n-1})$, 
the regularization $\bar{\cal J}(x)$ is evaluated as
\begin{equation}
\label{eq:barJn}
\bar{\cal J}(\bar\zeta)
=(1-q)^n\prod_{k=1}^n\frac{(q)_\8(qt^{-k})_\8\prod_{i=1}^2\prod_{j=1}^2(qa_i^{-1}b_j^{-1}t^{-(k-1)})_\8}
{(qt^{-1})_\8(qa_1^{-1}a_2^{-1}b_1^{-1}b_2^{-1}t^{-(n+k-2)})_\8}.
\end{equation}
\end{lem}
\noindent
{\bf Proof.} 
We denote by $C$ the right-hand side of (\ref{eq:barJn}). 
Then it is immediate to confirm that 
$C$ as a function of $a_j$ and $b_j$ satisfies the same $q$-difference equations as 
(\ref{eq:TaJ=cJ2}) and (\ref{eq:TbJ=cJ2}) of $\bar {\cal J}(\bar\zeta)$. Therefore 
the ratio $\bar {\cal J}(\bar\zeta)/C$ is invariant under the $q$-shift with respect to $a_j$ and $b_j$. 
\par
Next, for an integer $N$, let $T^N$ be the $q$-shift operator for a special direction defined as 
$$
T^N: b_1\to b_1q^{2N}, b_2\to b_2q^{-N}, a_1\to a_1q^{-N}, a_2\to a_2q^{-N}.
$$
Since we have 
\begin{equation*}
\bar\Phi(z)\Delta(z)=
\prod_{i=1}^n
z_i^{1-\alpha_1-\alpha_2-\beta_1-\beta_2-2(i-1)\tau}
\frac{(qb_1^{-1}z_i)_\8}{(a_1z_i)_\8}\frac{(qb_2^{-1}z_i)_\8}{(a_2z_i)_\8}
\prod_{1\le j<k\le n}
\frac{(qt^{-1}z_k/z_j)_\8}{(tz_k/z_j)_\8}
(1-z_k/z_j),
\end{equation*}
by definition $T^N\bar{J}(\bar\zeta)$ is written as
\begin{eqnarray*}
&&T^N\bar{J}(\bar\zeta)=(1-q)^n
\sum_{0\le \nu_1\le \nu_2\le \cdots\le \nu_n}\ 
\prod_{i=1}^n
(b_1t^{i-1}q^{\nu_i+2N})^{1-\alpha_1-\alpha_2-\beta_1-\beta_2-2(i-1)\tau+N}
\\[2pt]
&&\times
\frac{(qt^{i-1}q^{\nu_i})_\8}{(a_1b_1t^{i-1}q^{\nu_i+N})_\8}
\frac{(qb_2^{-1}b_1t^{i-1}q^{\nu_i+3N})_\8}{(a_2b_1t^{i-1}q^{\nu_i+N})_\8}
\prod_{1\le j<k\le n}
\frac{(qt^{-1+k-j}q^{\nu_k-\nu_j})_\8}{(t^{1+k-j}q^{\nu_k-\nu_j})_\8}
(1-t^{k-j}q^{\nu_k-\nu_j}),
\end{eqnarray*}
so that the leading term of the asymptotic behavior of $T^N\bar{J}(\bar\zeta)$ as $N\to +\8$ is given by
the term corresponding to $(\nu_1,\ldots,\nu_n)=(0,\ldots,0)$ in the above sum, which is
\begin{eqnarray}
\label{eq:TNJ}
&&\!\!\!\!\!\!
T^N\bar{J}(\bar\zeta)
\sim
(1-q)^n\prod_{i=1}^n
(b_1t^{i-1}q^{2N})^{1-\alpha_1-\alpha_2-\beta_1-\beta_2-2(i-1)\tau+N}
(qt^{i-1})_\8
\!\!\!
\prod_{1\le j<k\le n}
\!\!\!
\frac{(qt^{-1+k-j})_\8}{(t^{1+k-j})_\8}
\frac{(t^{k-j})_\8}{(qt^{k-j})_\8}\nonumber\\
&&\hskip 33pt
=(1-q)^n\prod_{i=1}^n
(b_1t^{i-1}q^{2N})^{1-\alpha_1-\alpha_2-\beta_1-\beta_2-2(i-1)\tau+N}
\frac{(q)_\8(t)_\8}{(t^{i})_\8}
\quad\quad(N\to +\8).
\end{eqnarray}
On the other hand, from (\ref{eq:barcalJ}), $\bar h(\bar\zeta)C$ is written as
\begin{eqnarray}
\label{eq:hC}
&&
\bar h(\bar\zeta)C=C\prod_{i=1}^n\frac{(b_1t^{i-1})^{1-\alpha_1-\alpha_2-\beta_1-\beta_2-2(i-1)\tau}}
{\theta(a_1b_1t^{i-1})\theta(a_2b_1t^{i-1})}
\frac{\theta(t)}{\theta(t^i)}\\
&&
=(1-q)^n
\prod_{i=1}^n\frac{(b_1t^{i-1})^{1-\alpha_1-\alpha_2-\beta_1-\beta_2-2(i-1)\tau}
(q)_\infty(t)_\infty(qa_1^{-1}b_2^{-1}t^{-(i-1)})_\8(qa_2^{-1}b_2^{-1}t^{-(i-1)})_\8}
{(t^i)_\infty(a_1b_1t^{i-1})_\infty(a_2b_1t^{i-1})_\infty
(qa_1^{-1}a_2^{-1}b_1^{-1}b_2^{-1}t^{-(n+i-2)})_\8},
\nonumber
\end{eqnarray}
so that we have 
\begin{eqnarray}
\label{eq:TNhC}
&&T^N\Big(\bar h(\bar\zeta)C\Big)
=(1-q)^n
\prod_{i=1}^n
(b_1t^{i-1}q^{2N})^{1-\alpha_1-\alpha_2-\beta_1-\beta_2-2(i-1)\tau+N}
\frac{(q)_\infty(t)_\infty}{(t^i)_\infty}
\nonumber\\
&&\hskip 40mm\times\frac{
(a_1^{-1}b_2^{-1}t^{-(i-1)}q^{1+2N})_\8(a_2^{-1}b_2^{-1}t^{-(i-1)}q^{1+2N})_\8}
{(a_1b_1t^{i-1}q^N)_\infty(a_2b_1t^{i-1}q^N)_\infty
(a_1^{-1}a_2^{-1}b_1^{-1}b_2^{-1}t^{-(n+i-2)}q^{1+N})_\8}
\nonumber\\[2pt]
&&
\sim (1-q)^n\prod_{i=1}^n
(b_1t^{i-1}q^{2N})^{1-\alpha_1-\alpha_2-\beta_1-\beta_2-2(i-1)\tau+N}
\frac{(q)_\infty(t)_\infty}{(t^i)_\infty}
\quad\quad(N\to +\8).
\end{eqnarray}
As we saw, the ratio $\bar {\cal J}(\bar\zeta)/C$ is invariant under the $q$-shift with respect to $a_j$ and $b_j$.
Thus $\bar {\cal J}(\bar\zeta)/C$ is also invariant under the $q$-shift $T^N$.
Therefore, comparing (\ref{eq:TNJ}) with (\ref{eq:TNhC}), we obtain
$$
\frac{\bar {\cal J}(\bar\zeta)}{C}
=T^N\frac{\bar {\cal J}(\bar\zeta)}{C}
=\frac{T^N\bar J(\bar\zeta)}{T^N \bar h(\bar\zeta)C}
=
\lim_{N\to +\infty}\frac{T^N\bar J(\bar\zeta)}{T^N \bar h(\bar\zeta)C}
=1,
$$
and thus $\bar {\cal J}(\bar\zeta)=C$. $\square$
\begin{cor}
\label{cor:barJzeta}
The truncated Jackson integral $\bar J(\bar\zeta)$ is evaluated as 
\begin{equation}\label{eq:barJzeta}
\begin{split}
\bar J(\bar\zeta)&=(1-q)^n
\prod_{i=1}^n\bigg[(b_1t^{i-1})^{1-\alpha_1-\alpha_2-\beta_1-\beta_2-2(i-1)\tau}\\
&\hskip 10pt\times
\frac{
(q)_\infty(t)_\infty(qa_1^{-1}b_2^{-1}t^{-(i-1)})_\8(qa_2^{-1}b_2^{-1}t^{-(i-1)})_\8}
{(t^i)_\infty(a_1b_1t^{i-1})_\infty(a_2b_1t^{i-1})_\infty
(qa_1^{-1}a_2^{-1}b_1^{-1}b_2^{-1}t^{-(n+i-2)})_\8}\bigg],
\end{split}
\end{equation}
\end{cor}
\noindent
{\bf Proof.} 
From Lemma \ref{lem:barJn}, $\bar J(\bar\zeta)=\bar h(\bar\zeta)\bar{\mathcal{J}}(\bar\zeta)$
is given by (\ref{eq:hC}). $\square$
\subsection{A remark on the relation between $\bar J(\bar\zeta)$ and $\bar I(\bar\zeta)$ of Aomoto's setting }
As an application of the $q$-difference equations (\ref{eq:TaJ=cJ}) and (\ref{eq:TbJ=cJ}) for $\bar J(x)$,  
we can show that 
the product formula (\ref{eq:01I(b)}) for $\bar I(\bar\zeta)$ in Corollary \ref{cor:bar I(a-N;b)}
(or (\ref{eq:I(zeta)2}) of $I(\zeta)$ in Proposition \ref{thm:aomoto} by the duality (\ref{eq:transform}) of parameters) is 
a special case of (\ref{eq:barJn}) in Lemma \ref{lem:barJn}. 
This indicates a way to prove the summation formula (\ref{eq:I(zeta)2}) from 
the product formula of the Jackson integral of Selberg type. 
\begin{cor}
\label{cor:04} For the point $x=\bar\zeta=(b_1,b_1t,\ldots,b_1t^{n-1})$, the truncated Jackson integral $\bar J(x)$ of 
Selberg type 
is expressed as 
\begin{equation}
\label{eq:02JandI}
\bar J(\bar\zeta)=
\bar I(\bar\zeta)
\prod_{i=1}^n\frac{(qa_1^{-1}b_{2}^{-1}t^{-(i-1)})_\8}{(b_1a_{2}t^{i-1})_\8},
\end{equation}
where $\bar I(\bar\zeta)$ is the truncated Jackson integral defined by {\rm (\ref{eq:01bar I(x)})} with the setting 
$\alpha=\alpha_{2}+\beta_{2}$. 
In particular, $\bar I(\bar\zeta)$ is expressed as {\rm (\ref{eq:01I(b)})}. 
\end{cor}
\begin{remark}
From (\ref{eq:02JandI}), $\bar I(\bar\zeta)$ is 
a limiting case of $\bar J(\bar\zeta)$ with the $q$-shift $a_{2}\to q^N a_{2}$
and $b_{2}\to q^{-N} b_{2}$ $(N\to +\8)$. 
Conversely, the product formula 
(\ref{eq:barJzeta}) of $\bar J(\bar\zeta)$ in 
Corollary \ref{cor:barJzeta} is reconstructed from 
the product formula (\ref{eq:01I(b)}) of $\bar I(\bar\zeta)$ 
via the connection (\ref{eq:02JandI}). 
\end{remark}
\noindent
{\bf Proof.} 
From (\ref{eq:TaJ=cJ}) and (\ref{eq:TbJ=cJ}) the recurrence relation of $\bar J(\bar\zeta)$ 
with respect to the $q$-shift $a_{2}\to qa_{2}$
and $b_{2}\to q^{-1} b_{2}$ is written as 
\begin{equation*}
\bar J(\bar\zeta)=T_{b_{2}}^{-1}T_{a_{2}}\bar J(\bar\zeta)\times 
\prod_{i=1}^{n}\frac{1-qa_1^{-1}b_{2}^{-1}t^{-(i-1)}}{1-b_1a_{2}t^{i-1}}.
\end{equation*}
By repeated use of this equation we have 
\begin{eqnarray}
\label{eq:02JandI1}
\bar J(\bar\zeta)=
T_{b_{2}}^{-N}T_{a_{2}}^N\bar J(\bar\zeta)
\prod_{i=1}^n\frac{(qa_1^{-1}b_{2}^{-1}t^{-(i-1)})_N}{(b_1a_{2}t^{i-1})_N}
=\lim_{N\to \8}T_{b_{2}}^{-N}T_{a_{2}}^N\bar J(\bar\zeta)
\prod_{i=1}^n\frac{(qa_1^{-1}b_{2}^{-1}t^{-(i-1)})_\8}{(b_1a_{2}t^{i-1})_\8}.
\end{eqnarray}
Moreover, by definition $\displaystyle \lim_{N\to \8}T_{b_{2}}^{-N}T_{a_{2}}^N\bar J(\bar\zeta)$ 
is written as 
\begin{eqnarray}
&&\lim_{N\to \8}T_{b_{2}}^{-N}T_{a_{2}}^N\bar J(\bar\zeta)
=\lim_{N\to \8}
(1-q)^n
\sum_{0\le\nu_1\le \nu_2\le \cdots \le \nu_n}
\prod_{i=1}^n
(b_1 t^{i-1} q^{\nu_i})^{1-\alpha_1-\alpha_2-\beta_1-\beta_2-2(i-1)\tau}
\nonumber\\
&&
\qquad\times\prod_{i=1}^n
\frac{(t^{i-1}q^{1+\nu_i})_\8(b_2^{-1}b_1t^{i-1}q^{1+\nu_i+N})_\8}
{(a_1b_1t^{i-1}q^{\nu_i})_\8(a_2b_1t^{i-1}q^{\nu_i+N})_\8}
\prod_{1\le j<k\le n}
\frac{(t^{k-j-1}q^{1+\nu_k-\nu_j})_\8}{(t^{k-j+1}q^{\nu_k-\nu_j})_\8}
(1-t^{k-j}q^{\nu_k-\nu_j})
\nonumber\\
&&=
(1-q)^n
\sum_{0\le\nu_1\le \nu_2\le \cdots \le \nu_n}
\prod_{i=1}^n
(b_1 t^{i-1} q^{\nu_i})^{1-\alpha_1-\alpha_2-\beta_1-\beta_2-2(i-1)\tau}
\frac{(t^{i-1}q^{1+\nu_i})_\8}
{(a_1b_1t^{i-1}q^{\nu_i})_\8}
\nonumber\\
&&\hskip 100pt
\times\prod_{1\le j<k\le n}
\frac{(t^{k-j-1}q^{1+\nu_k-\nu_j})_\8}{(t^{k-j+1}q^{\nu_k-\nu_j})_\8}
(1-t^{k-j}q^{\nu_k-\nu_j}),
\label{eq:02JandI2}
\end{eqnarray}
which exactly coincides with $\bar I(\bar\zeta)$ 
under the setting $\alpha=\alpha_2+\beta_2$. 
From (\ref{eq:02JandI1}) and (\ref{eq:02JandI2}), we therefore obtain (\ref{eq:02JandI}).

Next, using (\ref{eq:02JandI}) and (\ref{eq:barJzeta}) of Corollary \ref{cor:barJzeta},
the sum $
\bar I(\bar\zeta)
$ is conversely calculated as 
\begin{equation*}
\begin{split}
\bar I(\bar\zeta)
&=
\bar J(\bar\zeta)\prod_{i=1}^n\frac{(b_1a_2t^{i-1})_\8}{(qa_1^{-1}b_2^{-1}t^{-(i-1)})_\8}
\\
&=(1-q)^n
\prod_{i=1}^n\frac{
(b_1t^{i-1})^{1-\alpha_1-\alpha_2-\beta_1-\beta_2-2(i-1)\tau}
(q)_\8(t)_\8(qa_2^{-1}b_2^{-1}t^{-(i-1)})_\8}
{(t^{i})_\8(a_1b_1t^{i-1})_\8(qa_1^{-1}a_2^{-1}b_1^{-1}b_2^{-1}t^{-(n+i-2)})_\8},
\end{split}
\end{equation*}
which 
coincides with (\ref{eq:01I(b)}) of 
Corollary \ref{cor:bar I(a-N;b)} under the setting $\alpha=\alpha_2+\beta_2$.
$\square$
\appendix
\section{Appendix -- Derivation of the difference equations}
The aim of this section is to show a way to derive 
the $q$-difference equations (\ref{eq:TaJ=cJ}), (\ref{eq:TbJ=cJ}) in Proposition \ref{prop:q-diff01} 
for $\bar J(\zeta_i(b_1,b_2))$
using the shifted symmetric polynomials (the interpolation polynomials), which are defined as follows.
\begin{lem}[Knop--Sahi]
For $a,t\in \mathbb{C}^*$, $z=(z_1,z_2,\ldots,z_n)\in (\mathbb{C}^*)^n$, 
let $E_i(a;t;z)$ be the polynomials defined by
$$
E_r(a;t;z):=\sum_{1\le i_1< \cdots< i_r\le n}\prod_{k=1}^r (z_{i_k}-at^{i_k-k}) \quad\mbox{for}\quad r=1,2,\ldots, n,
$$
which are symmetric with respect to $z$, where $E_0(a;t;z)=1$. Then the polynomials $E_i(a;t;z)$ satisfy 
the vanishing property 
\begin{equation}
\label{eq:E-vanishing1}
E_i(a;t;\zeta_j)=0\quad\mbox{if}\quad 0\le j<i\le n,
\end{equation}
where
\begin{equation}
\label{eq:zeta-i}
\zeta_j:=(\hskip 1pt
z_1,z_2,\ldots,z_j,
\underbrace{a,at,\ldots,at^{n-j-1}}_{n-j}\hskip 1pt)\in (\mathbb{C}^*)^n.
\end{equation}
\end{lem}
\noindent
{\bf Proof.} See \cite[p.476, Proposition 3.1]{KS96}. $\square$ 
\begin{remark} 
The polynomials $E_i(a;t;z)$, $i=0,1,\ldots, n$, are called the {\em shifted symmetric polynomials} 
in the context \cite{KS96}. 
Using the factor theorem for the vanishing property (\ref{eq:E-vanishing1}) 
in this lemma, 
we immediately have 
\begin{equation}
\label{eq:Ei(zi)}
E_i(a;t;\zeta_i)=\prod_{k=1}^i(z_k-at^{n-i}).
\end{equation}
In particular, this is consistent with the $n$th symmetric polynomial
having the explicit form
\begin{equation}
E_n(a;t;z)=\prod_{k=1}^n (z_k-a). 
\end{equation}
\end{remark}

We rewrite the vanishing property (\ref{eq:E-vanishing1}) appropriately for the succeeding arguments of this section. 
\begin{lem}
\label{lem:vanishing}
Suppose that variables $z_1,\ldots,z_j$ in $\zeta_j$ are real numbers and satisfy 
\begin{equation}
\label{eq>>>}
z_1\gg z_2\gg \cdots\gg z_j\gg 0,
\end{equation}
which means $z_1/z_2\to \infty, z_2/z_3\to \infty, \ldots, z_{j-1}/z_j\to \infty$ and $z_j\to \infty$. 
Then the following asymptotic behavior holds: 
\begin{equation}
\label{eq:E-vanishing2}
\frac{E_i(a;t;z)\Delta(z)}{z_1^nz_2^{n-1}\cdots z_j^{n-j+1}}\Big|_{z=\zeta_j}
\sim\delta_{ij}\Delta^{\!\!(n-i)}(a,at,\ldots,at^{n-i-1}),
\end{equation}
where $\Delta^{\!\!(k)}(z_1,\ldots,z_k)$ denotes 
the difference product $\Delta(z)$ of $k$ variables.
\end{lem}
\noindent
{\bf Proof.} From (\ref{eq:E-vanishing1}), if $i>j$, 
then the left-hand side of (\ref{eq:E-vanishing2}) is exactly equal to 0. 
On the other hand, if $i<j$, then 
the degree of $E_i(a;t;z)\Delta(z)$ is lower than $z_1^nz_2^{n-1}\cdots z_j^{n-j+1}$, 
so that the left-hand side of (\ref{eq:E-vanishing2}) is estimated as $0$ under the condition (\ref{eq>>>}). 
If $i=j$, from (\ref{eq:Ei(zi)}), (\ref{eq:Delta}) and (\ref{eq:zeta-i}), we have 
$$
E_i(a;t;\zeta_i)\Delta(\zeta_i)=\Delta^{\!\!(n-i)}(a,at,\ldots,at^{n-i-1})\Delta^{\!\!(i)}(z_1,z_2,\ldots,z_i)
\prod_{j=1}^i\prod_{k=1}^{n-i+1}(z_j-at^{k-1}),
$$
which indicates (\ref{eq:E-vanishing2}). $\square$\\

We will state a key technical lemma for deriving $q$-difference equations. 
For this let $\bar\Phi(z)$ be the function defined by (\ref{eq:Phi2-1}) 
and for a function $\varphi(z)$, define the function $\nabla_{\!i}\varphi(z)$ ($1\le i\le n$) by
\begin{equation}
\label{eq:00nabla}
(\nabla_{\!i}\varphi)(z):=\varphi(z)-\frac{T_{z_i}\bar\Phi(z)}{\bar\Phi(z)}T_{z_i}\varphi(z),
\end{equation}
where $T_{z_i}$ means the shift operator of $z_i\to qz_i$, i.e.,
$T_{z_i}f(\ldots,z_i,\ldots)=f(\ldots,qz_i,\ldots)$. 
We then have  
\begin{lem}
\label{lem:00nabla=0}
For a meromorphic function $\varphi(z)$ on $({\mathbb C}^*)^n$, if the integral 
$$
\int_0^{\mbox{\small $x$}\8}\varphi(z)\bar\Phi(z)\varpi_q
$$
converges, then 
\begin{equation}
\label{eq:00nabla=0}
\int_0^{\mbox{\small $x$}\8}\bar\Phi(z)\nabla_{\!i}\varphi(z)\varpi_q=0. 
\end{equation}
Moreover, 
\begin{equation}
\label{eq:00A}
\int_0^{\mbox{\small $x$}\8}\bar\Phi(z){\cal A}\nabla_{\!i}\varphi(z)\varpi_q=0,
\end{equation}
where ${\cal A}$ indicates the skew-symmetrization defined in {\rm (\ref{eq:00Af})}. 
\end{lem}
\noindent
{\bf Proof.} 
From the definition (\ref{eq:00nabla}) of $\nabla_{\!i}$, 
(\ref{eq:00nabla=0}) is equivalent to the statement 
$$
\int_0^{\mbox{\small $x$}\8}\varphi(z)\bar\Phi(z)\varpi_q
=\int_0^{\mbox{\small $x$}\8}T_{z_i}\varphi(z)\, T_{z_i}\bar\Phi(z)
\varpi_q,
$$
if the left-hand side converges. And this equation 
is just confirmed from the fact that the Jackson integral is invariant under the $q$-shift $z_i\to qz_i$ ($1\le i\le n$).
Next we will confirm (\ref{eq:00A}). 
Taking account of the quasi-symmetry 
$
\sigma\bar\Phi(z)=U_\sigma(z)\bar\Phi(z)
$,
we have 
\begin{eqnarray*}
\bar\Phi(z){\cal A}\nabla_{\!i}\varphi(z)
&=&\bar\Phi(z)\sum_{\sigma\in S_n}(\sgn\, \sigma)\,\sigma(\nabla_{\!i}\varphi)(z)
=\sum_{\sigma\in S_n}(\sgn\, \sigma)U_\sigma(z)^{-1}\sigma\bar\Phi(z)\sigma(\nabla_{\!i}\varphi)(z)\\
&=&\sum_{\sigma\in S_n}(\sgn\, \sigma)U_\sigma(z)^{-1}\sigma\Big(\bar\Phi(z)\nabla_{\!i}\varphi(z)\Big).
\end{eqnarray*}
Since $U_\sigma(z)$ is invariant under the $q$-shift $z_i\to qz_i$ ($1\le i\le n$), we therefore obtain 
$$
\int_0^{\mbox{\small $x$}\8}\bar\Phi(z){\cal A}\nabla_{\!i}\varphi(z)\varpi_q
=\sum_{\sigma\in S}(\sgn\, \sigma)U_\sigma(x)^{-1}
\int_0^{\mbox{\small $x$}\8}\sigma\Big(\bar\Phi(z)\nabla_{\!i}\varphi(z)\Big)\varpi_q$$
$$
=\sum_{\sigma\in S}(\sgn\, \sigma)U_\sigma(x)^{-1}
\int_0^{\sigma^{-1}\mbox{\small $x$}\8}\!\!\!\!\!
\bar\Phi(z)\nabla_{\!i}\varphi(z)\varpi_q
=\sum_{\sigma\in S}(\sgn\, \sigma)U_\sigma(x)^{-1}\sigma\!\!
\int_0^{\mbox{\small $x$}\8}\bar\Phi(z)\nabla_{\!i}\varphi(z)\varpi_q,
$$
which vanishes from (\ref{eq:00nabla=0}). $\square$\\

We set 
\begin{equation}
\label{eq:ei=Ei}
e_i(a;t;z):=E_i(a;t;z^{-1}), 
\end{equation}
where $z^{-1}$ is specified by (\ref{eq:01x-1}). 
Since we have 
\begin{equation*}
\frac{T_{a_i}\bar\Phi(z)}{\bar\Phi(z)}
=\prod_{j=1}^n(z_j^{-1}-a_i)=e_n(a_i;t;z),
\quad 
\frac{T_{b_i}\bar\Phi(z)}{\bar\Phi(z)}
=\prod_{j=1}^n(z_j^{-1}-b_i^{-1})
=e_n(b_i^{-1};t^{-1};z),
\end{equation*}
the $q$-shifts of $\bar J(x)$ with respect to $a_i$ and $b_i$ are expressed by 
\begin{eqnarray}
T_{a_i}\bar J(x)
&=&\int_0^{\mbox{\small $x$}\8}e_n(a_i;t;z)\bar\Phi(z)\Delta(z)
\varpi_q,
\label{eq:Ta-barJ}
\\
T_{b_i}\bar J(x)
&=&\int_0^{\mbox{\small $x$}\8}e_n(b_i^{-1};t^{-1};z)\bar\Phi(z)\Delta(z)
\varpi_q.
\label{eq:Tb-barJ}
\end{eqnarray}
\begin{lem} 
\label{lem:Se=Se}
Suppose that $x=\zeta_k(b_1,b_2)$ $(k=0,1,\ldots, n)$, 
where $\zeta_k(x_1,x_2)$ is defined by {\rm (\ref{eq:zeta(x1x2)})}. 
Then the relation between $e_{i}(a_j;t;z)$ and $e_{i-1}(a_j;t;z)$ via the truncated Jackson integral is expressed as 
\begin{eqnarray}
\label{eq:Se=Se-a}
&&\int_0^{\mbox{\small $x$}} 
e_{i}(a_j;t;z)\bar\Phi(z)\Delta(z)\varpi_q\\
&&=
(-a_j)
\frac{(1-t^{n-i+1})(1-a_j^{-1}b_1^{-1}t^{-(n-i)})(1-a_j^{-1}b_2^{-1}t^{-(n-i)})}
{(1-t^{i})(1-a_1^{-1}a_2^{-1}b_1^{-1}b_2^{-1}t^{-(2n-i-1)})}
\int_0^{\mbox{\small $x$}} 
e_{i-1}(a_j;t;z)\bar\Phi(z)\Delta(z)\varpi_q, \quad
\nonumber
\end{eqnarray}
and the relation between $e_{i}(b_j^{-1};t^{-1};z)$ and $e_{i-1}(b_j^{-1};t^{-1};z)$ is expressed as
\begin{eqnarray}
\label{eq:Se=Se-b}
&&\int_0^{\mbox{\small $x$}} 
e_{i}(b_j^{-1};t^{-1};z)\bar\Phi(z)\Delta(z)\varpi_q\\
&&=
(-b_j^{-1})
\frac{(1-t^{-(n-i+1)})(1-a_1b_jt^{n-i})(1-a_2b_jt^{n-i})}
{(1-t^{-i})(1-a_1a_2b_1b_2t^{2n-i-1})}
\int_0^{\mbox{\small $x$}} 
e_{i-1}(b_j^{-1};t^{-1};z)\bar\Phi(z)\Delta(z)\varpi_q. \quad
\nonumber
\end{eqnarray}
\end{lem}
\begin{remark}
The relations (\ref{eq:Se=Se-a}) and (\ref{eq:Se=Se-b}) are identical upon the interchange of parameters as 
$(a_1,a_2,b_1,b_2,t)\to (b_1^{-1},b_2^{-1},a_1^{-1},a_2^{-1},t^{-1})$. 
\end{remark}
\begin{remark}
\label{remark:Se=Se}
By repeated use of (\ref{eq:Se=Se-a}), 
from (\ref{eq:Ta-barJ}), we immediately obtain the $q$-difference equation (\ref{eq:TaJ=cJ})  
presented in Proposition \ref{prop:q-diff01}. In the same manner,  
the $q$-difference equation (\ref{eq:TbJ=cJ}) in Proposition \ref{prop:q-diff01} is deduced from (\ref{eq:Se=Se-b})
using (\ref{eq:Tb-barJ}). 
\end{remark}

The rest of this subsection is devoted to the proof of the above lemma. 
We will show a further lemma before proving Lemma \ref{lem:Se=Se}.
For this purpose we abbreviate $e_i(a_1;t;z)$ and $E_i(a_1;t;z)$ by $e_i(z)$ and $E_i(z)$, respectively, 
and $(k)$ of $e_i^{(k)}(z)$, $E_i^{(k)}(z)$, $\Delta^{\!\!(k)}(z)$ means that these functions are of $k$ variables. 
We also use the symbol $(\widehat{z}_i):=(z_1,\ldots,z_{i-1},z_{i+1},\ldots,z_n)\in (\mathbb{C}^*)^{n-1}$ for $i=1,\ldots,n$.
\begin{lem} 
\label{lem:Anabla-phi}
Put 
\begin{equation}
\label{eq:phi-1}
\phi(z):=z_1^{-1}(z_1-b_1)(z_1-b_2)\prod_{k=2}^n(z_1-tz_k).
\end{equation}
Then
\begin{equation}
\label{eq:Anabla-phi-1}
(-1)^{n-1}{\cal A}\Big[\nabla_{\!1}\Big(\phi(z)\,e_{i-1}^{(n-1)}(\widehat{z}_1)\Delta^{\!\!(n-1)}(\widehat{z}_1)\Big)\Big]
=\Big(c_{i}e_{i}^{(n)}(z)+c_{i-1}e_{i-1}^{(n)}(z)\Big)\Delta^{\!\!(n)}(z),
\end{equation}
where the coefficients $c_{i-1}$ and $c_i$ are given by 
\begin{eqnarray}
c_{i-1}&=&(n-1)!(-1)^{n-1}a_1^{-1}t^{-(n-i)}(1-a_1b_1t^{n-i})(1-a_1b_2t^{n-i})(1-t^{n-i+1})/(1-t),
\qquad
\label{eq:coeffi-(i-1)}\\[2pt]
c_i&=&(n-1)!(-1)^{n}a_1^{-1}a_2^{-1}t^{-(n-1)}(1-a_1a_2b_1b_2t^{2n-i-1})(1-t^i)/(1-t).
\label{eq:coeffi-(i)}
\end{eqnarray}
\end{lem}
\noindent
{\bf Proof.} We initially compute $\nabla_{\!1}\phi(z)$. By definition $T_{z_1}\bar\Phi(z)/\bar\Phi(z)$ is written as 
\begin{eqnarray}
\label{eq:Tz1Phi/Phi}
\frac{T_{z_1}\bar\Phi(z)}{\bar\Phi(z)}
&=&(a_1a_2b_1b_2q^{n-2})^{-1}\frac{(1-a_1z_1)(1-a_2z_1)}{(1-qb_1^{-1}z_1)(1-qb_2^{-1}z_1)}
\prod_{i=2}^n \frac{1-t^{-1}z_i/z_1}{1-q^{-1}tz_i/z_1}
\nonumber\\
&=&\frac{q^{-n}z_1^{-n}(1-a_1^{-1}z_1^{-1})(1-a_2^{-1}z_1^{-1})}
{z_1^{-n}(1-q^{-1}b_1z_1^{-1})(1-q^{-1}b_2z_1^{-1})}
\prod_{k=2}^n\frac{1-t^{-1}z_1^{-1}/z_k^{-1}}{1-q^{-1}tz_1^{-1}/z_k^{-1}}.
\end{eqnarray}
Since the function $\phi(z)$ in (\ref{eq:phi-1}) is rewritten as 
\begin{equation}
\label{eq:phi-2}
\phi(z)=\frac{1}{z_1^{-n}}(1-b_1z_1^{-1})(1-b_2z_1^{-1})\prod_{k=2}^n(1-tz_1^{-1}/z_k^{-1}),
\end{equation}
using (\ref{eq:00nabla}) and (\ref{eq:Tz1Phi/Phi}) the function $\nabla_{\!1}\phi(z)$ is computed explicitly as 
\begin{equation}
\label{eq:nabra-phi}
\nabla_{\!1}\phi(z)=\phi(z)-
\frac{1}{z_1^{-n}}(1-a_1^{-1}z_1^{-1})(1-a_2^{-1}z_1^{-1})\prod_{k=2}^n(1-t^{-1}z_1^{-1}/z_k^{-1}). 
\end{equation}

Next, using (\ref{eq:ei=Ei}), the equation (\ref{eq:Anabla-phi-1}) is transformed by $z\to z^{-1}$, i.e., 
$$
(-1)^{n-1}{\cal A}\Big[\Big(\nabla_{\!1}\phi(z^{-1})\Big)
E_{i-1}^{(n-1)}(\widehat{z}_1)
\Delta^{\!\!(n-1)}(\widehat{z}_1^{\,-1})\Big]
=\Big(c_{i}E_{i}^{(n)}(z)+c_{i-1}E_{i-1}^{(n)}(z)\Big)\Delta^{\!\!(n)}(z^{-1}), 
$$
which is rewritten as 
$$
(-1)^{n-1}{\cal A}\bigg[\nabla_{\!1}\phi(z^{-1}) 
E_{i-1}^{(n-1)}(\widehat{z}_1)\frac{(-1)^{n-1\choose 2}\Delta^{\!\!(n-1)}(\widehat{z}_1)}{(z_2\cdots z_n)^{n-2}}\bigg]=
\Big(c_{i}E_{i}^{(n)}(z)+c_{i-1}E_{i-1}^{(n)}(z)\Big)\frac{(-1)^{n\choose 2}\Delta^{\!\!(n)}(z)}{(z_1z_2\cdots z_n)^{n-1}},
$$
so that it suffices to prove the following instead of (\ref{eq:Anabla-phi-1}):
\begin{equation}
\label{eq:Anabla-phi-2}
{\cal A}\Big[z_1^{n-1}z_2\cdots z_n\nabla_{\!1}\phi(z^{-1})\,E_{i-1}^{(n-1)}(\widehat{z}_1)\Delta^{\!\!(n-1)}(\widehat{z}_1)\Big]=
\Big(c_{i}E_{i}^{(n)}(z)+c_{i-1}E_{i-1}^{(n)}(z)\Big)\Delta^{\!\!(n)}(z). 
\end{equation}

We will prove the above equation. If we put $$
\bar\phi_i(z):=z_1^{n-1}z_2\cdots z_n\nabla_{\!1}\phi(z^{-1}) E_{i-1}^{(n-1)}(\widehat{z}_1)
\Delta^{\!\!(n-1)}(\widehat{z}_1),
$$ then, from (\ref{eq:nabra-phi}) $\bar\phi_i(z)$ is computed as 
\begin{eqnarray}
\label{eq:barphiiz}
&&\bar\phi_i(z)=
z_1^{-1}\bigg[
(1-b_1z_1)(1-b_2z_1)\prod_{k=2}^n(z_k-tz_1)\nonumber\\
&&\hskip 70pt 
-(1-a_1^{-1}z_1)(1-a_2^{-1}z_1)\prod_{k=2}^n(z_k-t^{-1}z_1)\bigg]
E_{i-1}^{(n-1)}(\widehat{z}_1)
\Delta^{\!\!(n-1)}(\widehat{z}_1),\qquad
\end{eqnarray}
which is a polynomial of $z_1,\ldots,z_n$. Taking account of the degree of the polynomial $\bar\phi_i(z)$, 
$${\cal A}(z_1z_2\cdots z_{i}\times z_1^{n-1}z_2^{n-2}\cdots z_{n-1})$$
is the term of highest degree in the skew-symmetrization ${\cal A}\bar\phi_i(z)$, 
which is thus expanded as 
\begin{equation}
\label{eq:Aphi-expand1}
{\cal A}\bar\phi_i(z)=\sum_{j=0}^{i} c_{j}E_j^{(n)}(z)\Delta^{\!\!(n)}(z),
\end{equation}
where $c_j$ are some constants. 
For any $x\in \mathbb{R}$ we set
$$
\xi_j:=(\hskip 1pt\underbrace{x^{j},x^{j-1},\ldots,x^2,x}_j,
\underbrace{a_1,a_1t,\ldots,a_1t^{n-j-1}}_{n-j}\hskip 1pt)\in (\mathbb{C}^*)^n,
$$
which is a special case of $\zeta_j$ specified by (\ref{eq:zeta-i}). 
Then, from Lemma \ref{lem:vanishing}, $c_j$ in (\ref{eq:Aphi-expand1}) is written as 
\begin{equation}
\label{eq:Aphi-expand1.5}
\lim_{x\to \infty}\frac{{\cal A}\bar\phi_i(z)}{z_1^nz_2^{n-1}\cdots z_j^{n-j+1}}\Big|_{z=\xi_j}
=c_j\lim_{x\to \infty}\frac{E_j^{(n)}(z)\Delta^{\!\!(n)}(z)}{z_1^nz_2^{n-1}\cdots z_j^{n-j+1}}\Big|_{z=\xi_j}
=c_j \,\Delta^{\!\!(n-j)}(a_1,a_1t,\ldots,a_1t^{n-j-1}).
\end{equation}

On the other hand, from the explicit form (\ref{eq:barphiiz}) of $\bar\phi_i(z)$, ${\cal A}\bar\phi_i(z)$ is also 
expressed as 
\begin{equation}
\label{eq:Aphi-expand2}
{\cal A}\bar\phi_i(z)=(n-1)!\sum_{j=1}^n (-1)^{j-1}\Big(F_j(z)-G_j(z)\Big)
\Delta^{\!\!(n-1)}(\widehat{z}_j)E_{i-1}^{(n-1)}(\widehat{z}_j),
\end{equation}
where
\begin{eqnarray}
F_j(z)&=&z_j^{-1}(1-b_1z_j)(1-b_2z_j)\prod_{1\le k\le n\atop k\ne j}(z_k-tz_j),
\label{eq:Fj(z)}\\
G_j(z)&=&z_j^{-1}(1-a_1^{-1}z_j)(1-a_2^{-1}z_j)\prod_{1\le k\le n\atop k\ne j}(z_k-t^{-1}z_j),
\label{eq:Gj(z)}
\end{eqnarray}
which satisfy the vanishing property 
\begin{equation}
\label{eq:FGvanishing}
\begin{array}{c}
F_j(\xi_k)=0\quad\mbox{if}\quad k< j< n, \\[5pt]
G_j(\xi_k)=0\quad\mbox{if}\quad k< j\le n,
\end{array}
\end{equation}
and the evaluation
\begin{equation}
\label{eq:Fn-evaluation}
\lim_{x\to \infty}
\frac{F_{n}(z)}
{z_1z_2\cdots z_k}\Big|_{z=\xi_k}
=\frac{(1-a_1b_1t^{n-k-1})(1-a_1b_2t^{n-k-1})(1-t^{n-k})}{a_1t^{n-k-1}(1-t)}
\prod_{j=1}^{n-k-1}(a_1t^{j-1}-a_1t^{n-k-1}).
\end{equation}
Notice that, from Lemma \ref{lem:vanishing}, 
$E_{i-1}^{(n-1)}(\widehat{z}_j)\Delta^{\!\!(n-1)}(\widehat{z}_j)$ 
satisfies the vanishing properties
\begin{equation}
\label{eq:E-vanishing-2}
\lim_{x\to \infty}
\frac{E_{i-1}^{(n-1)}(\widehat{z}_j)\Delta^{\!\!(n-1)}(\widehat{z}_j)}{z_1^{n-1}\cdots z_{j-1}^{n-j+1}
z_{j+1}^{n-j}\cdots z_k^{n-k+1}}\Big|_{z=\xi_k}
=\delta_{ik}\Delta^{\!\!(n-k)}(a_1,a_1t,\ldots,a_1t^{n-k-1})
\quad\mbox{if}\quad 1\le j\le k,
\end{equation}
and
\begin{equation}
\label{eq:E-vanishing-3}
\lim_{x\to \infty}
\frac{E_{i-1}^{(n-1)}(\widehat{z}_n)\Delta^{\!\!(n-1)}(\widehat{z}_n)}
{z_1^{n-1}z_2^{n-2}\cdots z_{k}^{n-k}}
\Big|_{z=\xi_k}
=\delta_{i-1,k}\Delta^{\!\!(n-k-1)}(a_1,a_1t,\ldots,a_1t^{n-k-2}).
\end{equation}

We now prove the expression (\ref{eq:Anabla-phi-2}). 
If $0\le k\le i-2$, then, from (\ref{eq:FGvanishing}), (\ref{eq:E-vanishing-2}) and (\ref{eq:E-vanishing-3}) 
the equation (\ref{eq:Aphi-expand2}) indicates 
\begin{equation}
\label{eq:Aphi-expand3}
\lim_{x\to \infty}\frac{{\cal A}\bar\phi_i(z)}{z_1^nz_2^{n-1}\cdots z_k^{n-k+1}}\Big|_{z=\xi_k}
=0.
\end{equation}
Comparing (\ref{eq:Aphi-expand1.5}) with (\ref{eq:Aphi-expand3}), we obtain $c_0=c_1=\cdots=c_{i-2}=0$, 
which means (\ref{eq:Anabla-phi-2}) holds. 

Next we will evaluate $c_{i-1}$. 
From (\ref{eq:FGvanishing}), (\ref{eq:Fn-evaluation}), (\ref{eq:E-vanishing-2}) and (\ref{eq:E-vanishing-3}), 
\begin{eqnarray}
\label{eq:Aphi-expand4}
&&\lim_{x\to \infty}\frac{{\cal A}\bar\phi_i(z)}{z_1^nz_2^{n-1}\cdots z_{i-1}^{n-i+2}}\Big|_{z=\xi_{i-1}}
=(n-1)!(-1)^{n-1}
\lim_{x\to \infty}
\frac{F_{n}(z)}
{z_1z_2\cdots z_{i-1}}
\frac{E_{i-1}^{(n-1)}(\widehat{z}_n)\Delta^{\!\!(n-1)}(\widehat{z}_n)}
{z_1^{n-1}z_2^{n-2}\cdots z_{i-1}^{n-i+1}}
\Big|_{z=\xi_{i-1}}
\nonumber\\
&&
\qquad=(n-1)!(-1)^{n-1}
\frac{(1-a_1b_1t^{n-i})(1-a_1b_2t^{n-i})(1-t^{n-i+1})}
{a_1t^{n-i}(1-t)}
\prod_{j=1}^{n-i}(a_1t^{j-1}-a_1t^{n-i})
\nonumber\\[-3pt]
&&\qquad\quad
\times
\Delta^{\!\!(n-i)}(a_1,a_1t,\ldots,a_1t^{n-i-1}).
\end{eqnarray}
Comparing (\ref{eq:Aphi-expand1.5}) with (\ref{eq:Aphi-expand4}) using the relation
$$
\frac{\Delta^{\!\!(n-i+1)}(a_1,a_1t,\ldots,a_1t^{n-i})}{\Delta^{\!\!(n-i)}(a_1,a_1t,\ldots,a_1t^{n-i-1})}
=\prod_{j=1}^{n-i}(a_1t^{j-1}-a_1t^{n-i}),
$$
we therefore obtain the explicit expression of $c_{i-1}$ as (\ref{eq:coeffi-(i-1)}).

Lastly we will evaluate $c_{i}$. 
From the explicit forms (\ref{eq:Fj(z)}) and (\ref{eq:Gj(z)}) of $F_j(z)$ and $G_j(z)$, we have 
$$
\lim_{x\to \infty}\frac{F_j(z)-G_j(z)}{z_1z_2\cdots z_{j-1}z_j^{n-j+1}}\Big|_{z=\xi_i}
=(-1)^{n-j}(b_1b_2t^{n-j}-a_1^{-1}a_2^{-1}t^{-(n-j)})
\quad\mbox{if}\quad 1\le j\le i. 
$$
Using (\ref{eq:FGvanishing}), (\ref{eq:E-vanishing-2}), (\ref{eq:E-vanishing-3}) and the above evaluation, we have 
\begin{eqnarray}
\label{eq:Aphi-expand5}
&&\lim_{x\to \infty}\frac{{\cal A}\bar\phi_i(z)}{z_1^nz_2^{n-1}\cdots z_i^{n-i+1}}\Big|_{z=\xi_i}
\nonumber\\
&&=
(n-1)!\sum_{j=1}^i(-1)^{j-1}\lim_{x\to \infty}\frac{F_j(z)-G_j(z)}{z_1z_2\cdots z_{j-1}z_j^{n-j+1}}
\frac{E_{i-1}^{(n-1)}(\widehat{z}_j)\Delta^{\!\!(n-1)}(\widehat{z}_j)}{z_1^{n-1}\cdots z_{j-1}^{n-j+1}
z_{j+1}^{n-j}\cdots z_i^{n-i+1}}\Big|_{z=\xi_i}
\nonumber\\
&&=(n-1)!(-1)^{n-1}\Delta^{\!\!(n-i)}(a_1,a_1t,\ldots,a_1t^{n-i-1})\sum_{j=1}^i(b_1b_2t^{n-j}-a_1^{-1}a_2^{-1}t^{-(n-j)})\nonumber\\
&&=(n-1)!(-1)^{n}\frac{}{}\frac{(1-a_1a_2b_1b_2t^{2n-i-1})(1-t^i)}{a_1a_2t^{n-1}(1-t)}\Delta^{\!\!(n-i)}(a_1,a_1t,\ldots,a_1t^{n-i-1}). 
\end{eqnarray}
Comparing (\ref{eq:Aphi-expand1.5}) with (\ref{eq:Aphi-expand5}), 
we therefore obtain the explicit expression of $c_{i}$ as (\ref{eq:coeffi-(i)}). $\square$\\

\noindent
{\bf Proof of Lemma \ref{lem:Se=Se}.} 
We will prove (\ref{eq:Se=Se-a}) for $a_j$ first. 
Without loss of generality, it suffices to show (\ref{eq:Se=Se-a}) for $a_1$. 
Suppose that $x=\zeta_k(b_1,b_2)$, $k=0,1,\ldots, n$. 
If we set 
$\varphi(z)=\phi(z)\,e_{i-1}^{(n-1)}(\widehat{z}_1)\Delta^{\!\!(n-1)}(\widehat{z}_1)$, 
where $\phi(z)$ is defined by (\ref{eq:phi-1}), then 
the truncated Jackson integral 
$
\int_0^{\mbox{\small $x$}}\varphi(z){\bar\Phi}(z)\varpi_q
$
converges absolutely if $|a_1a_2b_1b_2|$ is sufficiently large. 
Therefore, applying (\ref{eq:00A}) in Lemma \ref{lem:00nabla=0} to 
the fact (\ref{eq:Anabla-phi-1}) in Lemma \ref{lem:Anabla-phi}, we obtain the relation
$c_i\int_0^{\mbox{\small $x$}} 
e_{i}(a_1;t;z)\bar\Phi(z)\Delta(z)\varpi_q
+c_{i-1}
\int_0^{\mbox{\small $x$}} 
e_{i-1}(a_1;t;z)\bar\Phi(z)\Delta(z)\varpi_q=0,$ 
where $c_{i-1}$ and $c_i$ are given in (\ref{eq:coeffi-(i-1)}) and (\ref{eq:coeffi-(i)}), respectively. 
This relation coincides with (\ref{eq:Se=Se-a}). 

Next we will show (\ref{eq:Se=Se-b}) for $b_j$ of the case $j=1$ 
in the same manner as above. 
Here if we exchange $a_1,a_2$ and $t$ with $b_1^{-1},b_2^{-1}$ and $t^{-1}$, respectively,   
in the above proof of (\ref{eq:Se=Se-a})
including that of Lemma \ref{lem:Anabla-phi}, 
the way of argument is completely symmetric for this exchange. 
Therefore (\ref{eq:Se=Se-b}) for $b_1$ is obtained 
exchanging $a_1,a_2$ and $t$ with $b_1^{-1},b_2^{-1}$ and $t^{-1}$, 
respectively, on the coefficient of  (\ref{eq:Se=Se-a}). 
$\square$

\begin{remark} In the above proof, the assumption $x=\zeta_k(b_1,b_2)$, $k=0,1,\ldots, n$, 
for the truncated Jackson integral 
$
\int_0^{\mbox{\small $x$}}\varphi(z){\bar\Phi}(z)\varpi_q
$
is necessary from the technical view point. 
In the case for any $x\ne\zeta_i(b_1,b_2)$, taking account of the influence of the terms 
$\varphi(xq^\nu){\bar\Phi}(xq^\nu)$, $\nu\not\in \Lambda_i$, the convergence of 
$
\int_0^{\mbox{\small $x$}\infty}\varphi(z){\bar\Phi}(z)\varpi_q
$
is very subtle, and generally it is not assured. 
Lemma \ref{lem:00nabla=0} requires this convergence. 
\end{remark}

\subsection*{Acknowledgements}
We thank Prof.~J.V.~Stokman for altering us to the works \cite{St00, TV97} in relation
to our Corollary \ref{cor:main1}.
This work was supported by the Australian Research Council (Grant DP110102317) and JSPS KAKENHI Grant Number 25400118.

\end{document}